\documentclass[%prl,
twocolumn,
nofootinbib,
 amsmath,amssymb,
 aps,
pra,
]{revtex4-2}

\usepackage{graphicx} % Required for inserting images
\usepackage{amsmath}
\usepackage{amsfonts}
\usepackage{amssymb}
\usepackage{color}
\usepackage{bm}
\usepackage{fullpage}
\usepackage{comment}
\usepackage{todonotes}
\usepackage{tikz-cd} 
\usetikzlibrary{arrows,arrows.meta}

\begin{document}

\newcommand{\mycorr}[1]{ \textcolor{black}{#1}} 
\title{Learning Global Linear Representations of Nonlinear Dynamics}

\author{Thomas Breunung\textsuperscript{\textdagger}}
\email{breunung@wisc.edu}
\affiliation{Department of Mechanical Engineering, University of Maryland, College Park, College Park, Maryland, USA}
\thanks{Current address: Department of Mechanical Engineering, University of Wisconsin - Madison, Madison, Wisconsin, USA}

\author{Florian Kogelbauer}
\email{floriank@ethz.ch} 
\affiliation{Department of Mechanical and Process Engineering, ETH Zurich, Switzerland}

\date{\today}

%\begin{abstract}
%\mycorr{We address the question of global linearization of nonlinear dynamical systems. Using neural networks, we demonstrate that nonlinear phenomena can be represented by low-dimensional immersions using special linear reference dynamics. We illustrate our method on several examples, including nested families of periodic orbits, stable limit cycles and the coexistence of multiple steady states in forced-damped oscillators and draw connections to Koopman theory. The obtained linearizations are finite dimensional and exceed the phase space dimension of the underlying linear system by one at most.  }

\begin{abstract}
\mycorr{While linear systems are well-understood, no explicit solution for general nonlinear systems exists. A classical approach to make the understanding of linear system available in the nonlinear setting is to represent a nonlinear system by a linear model. While progress has been made in extending linearization techniques to larger domains and more complex attractor geometries, recent work has highlighted the limitations of these techniques when applied to nonlinear dynamics, such as those with coexisting attractors. In this work, we show nonlinear dynamics with a continuous Koopman spectrum, a limit cycle, and coexisting solutions that can be globally linearized. To this end, we explicitly construct linear systems mimicking these nonlinear behaviors. Subsequently, we approximate transformations between linear and nonlinear systems with deep neural networks. This approach yields finite dimensional linearizations exceeding the phase space dimension of the underlying linear system by one at most. }

\smallbreak
\end{abstract}

\maketitle

\section{Introduction}

Linear dynamical systems constitute the most basic and well-understood class of \mycorr{differential equations}. Linearity enables powerful principles, such as superposition and the derivation of closed-form solutions \cite{arnol2013mathematical}. Many well-established and computationally efficient algorithms for control, estimation, and prediction are thus based on linearizations \cite{guckenheimer2013nonlinear}. On the other hand, many phenomena in nature and engineering are inherently nonlinear, from shock formation \cite{whitham2011linear} to strange attractors \cite{lorenz1963deterministic} and turbulence in fluids \cite{holmes1997low}. While classical results on the relation between linearized and nonlinear dynamics, such as the center manifold theorem, are local in nature \cite{sijbrand1985properties}, a global description of nonlinear systems by linear principles would be highly desirable. This motivates our fundamental question: Can nonlinear systems be accurately and globally represented by linear dynamics? 

Recently, this basic question has drawn considerable attention in the context of linear immersions and lies at the heart of Koopman theory \cite{brunton2021modern}. In his seminal work \cite{koopman1931hamiltonian}, Koopman formulated classical Hamiltonian mechanics as unitary transformations in Hilbert spaces by representing  nonlinear flow maps as the action of linear operators. This result has then been extended to general finite-dimensional, nonlinear systems \cite{brunton2021modern}. 

In Koopman theory, dynamical features of the flow map are encoded in the spectrum of a transition operator, called Koopman operator. The corresponding generalized eigenfunctions constitute fundamental buildings blocks of the nonlinear dynamics \cite{budivsic2012geometry}. Indeed, once computed, the eigenfunctions of the Koopman operator provide an intrinsic coordinate system that allows for a uniform linearization of nonlinear dynamics \cite{mezic2020spectrum} on their domain of definition and the evaluation of the flow map is almost instantaneous. Dynamic Mode Decomposition (DMD)~\cite{schmid2010dynamic} and its nonlinear extensions~\cite{williams2016extending} serve as a mechanism that can provide a finite-dimensional approximation of Koopman eigenfunctions.

In the context of operator dynamics, the adjoint of the Koopman operator, called \textit{transfer operator} or \textit{ 
Perron–Frobenius operator}~\cite{ruelle2004thermodynamic}, is also frequently used in data-driven dimensional reduction of dynamical systems \cite{klus2018data}. Particular relevance is attributed to its eigenstructure - especially the existence of a leading eigenvector with strictly positive eigenvalue - implying the existence of an invariant measure of the underlying dynamical system \cite{baladi2000positive}. 

Recently, the application of machine learning to Koopman theory has gained considerable attention \cite{bevanda2021koopman,korda2020optimal}. Although comparably young as a field, there already exists a vast literature on the application of neural networks and related techniques to Koopman theory. In their influential work \cite{lusch2018deep}, Lusch \textit{et al.} discuss the application of deep learning techniques to Koopman theory. Further research directions include data-driven discovery of Koopman eigenfunctions for control problems \cite{kaiser2021data}, nonlinear normal modes in forced system \cite{rostamijavanani2023study} and data driven construction of Koopman eigenfunctions \cite{korda2020optimal} - only to name a few. 

One of the main difficulties of Koopman theory and related techniques, however, lies in the computation of eigenfunctions \cite{williams2014kernel}. Even for comparably simple systems, Koopman eigenfuctions often remain uninterpretable and complex~\cite{lusch2018deep}. Furthermore, the choice of basis functions to approximate the infinite-dimensional Koopman operator as a large matrix is by no means trivial and may limit the representation of possible dynamics considerably \cite{yeung2019learning}. Physical insight in the dynamics might facilitate the search of appropriate observables to pick as basis functions, while the selection process might be automated \cite{brunton2016koopman}. However, it remains, in principle, unclear, if there are theoretically guaranteed justifications for linear representations of nonlinear dynamics in the Koopman framework all together. 

From a more abstract point of view, the question of linear representation of nonlinear flows has received wide attention in the theory of dynamical systems. These results can serve as an underpinning to justify computations of Koopman eigenfunctions. For example, the classical Hartman--Grobman theorem \cite{katok1995introduction} proves topological conjugacy of a nonlinear flow map to its linearization around a hyperbolic fixed point. Various extensions of this theorem to normally hyperbolic or attracting invariant manifolds exist~\cite{eldering2018global}. In addition, assuming certain non-resonance conditions, the Sternberg linearization theorem guarantees smooth local equivalence of a non-linear system to its linearization around a hyperbolic fixed-point~\cite{sternberg1958structure}. In this setting one may employ the method of normal forms~\cite{murdock2003normal}, wherein through series of smooth transformations
one obtains a set of coordinates in which the original dynamics are linear\footnote{We note that the method of normal forms can also yield nonlinear terms in the essential dynamics, i.e. the normal form. However, if the Sternberg linearization theorem applies the associated normal form is linear.}. Rather recently it has been shown that the dynamics of a stable fixed point can be linearized within this fixed point's basin of attraction \cite{lan2013linearization}. This result has been used to construct two-dimensional limit cycles (and quasi-periodic solutions) in autonomous systems \cite{mezic2020spectrum}. In this context, we also mention Carleman linearization \cite{kowalski1991nonlinear}, which is based on the Taylor expansion of the right-hand side of nonlinear dynamics. We emphasize that these linearization results differ from Koopman theory in two fundamental aspects. Firstly, those results establish a  topological or differential conjugacy, i.e., linear representation in the same dimension as the nonlinear dynamics, while Koopman linearizations are generally infinite-dimensional even for finite dimensional systems~\cite{brunton2021modern}. Moreover, the aforementioned special coordinate systems and linearization are, in general, only guaranteed to exist locally, i.e. in a neighborhood of fixed points, (quasi-) periodic solutions, or normally-hyperbolic attracting sets \cite{mezic2020spectrum}.

A major obstruction to establishing global linear conjugacy consists in nonlinear phenomena, such as continuous \mycorr{Koopman} spectra, limit cycles, or the coexistence of multiple steady states. Indeed, Liu \textit{et al.} \cite{liu2023properties} show the limitations of continuous immersions of systems with multiple limit sets. 
\mycorr{It is well-known \cite{smale1967differentiable} that dynamical systems with coexisting steady states are non-linearizable by a continuous homeomorphism. This potential limitation to global linear immersion was also pointed out more recently in the context of Koopman theory \cite{cenedese2022data,kvalheim2023linearizability}.}
Furthermore, Page and Kershwell \cite{page2019koopman} indicate that dynamical systems with multiple invariant solutions (such as Couette flows or systems with heteroclinic connections) cannot be represented by globally convergent Koopman expansions. It is thus by no means obvious whether linearizing transformations can capture nonlinear phenomena.

\begin{figure}[h!]
 	\begin{center}
			\includegraphics[width=0.5\textwidth]{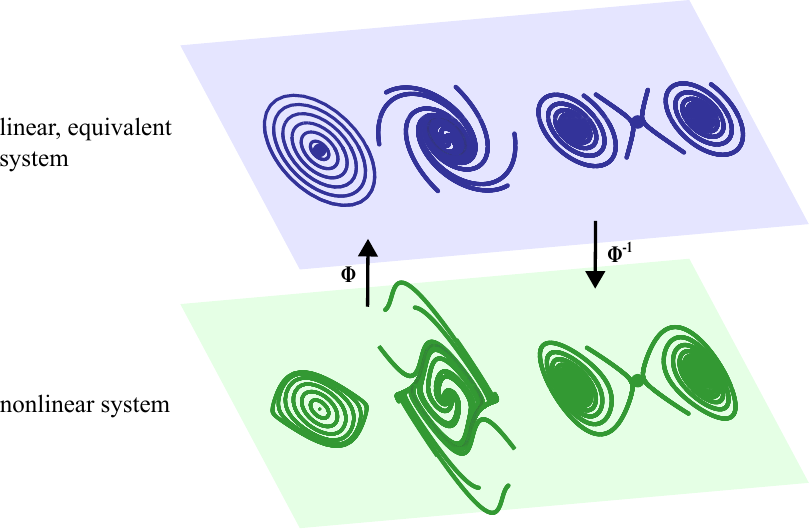}
 	\end{center}
 	\caption{Illustration of our approach: We learn a mapping $\bm{\Phi}$ (together with its restricted inverse $\bm{\Phi}^{-1}$) from trajectory data that globally embeds a nonlinear system into a linear system with only slightly larger dimension. Specifically, we consider systems with continuous \mycorr{Koopman} spectra, limit cycles, and coexisting steady state solutions and provide prototypical conjugated linear systems.
   }
 	\label{fig:overview}
\end{figure} 

We are thus prompted with the question whether nonlinear phenomena can be accurately captured by global linear immersions (as opposed to local linear conjugacy). In addition, we seek to clarify if these linearization spaces are necessarily infinite dimensional, as indicated in Koopman theory, or if low dimensional linearizations are achievable. 

\mycorr{In this work, we show examples of nonlinear systems that can be globally linearized by explicit constructions and by learning global linearizations from trajectory data.} Our approach is based on prototypical linear dynamics that mimic specific nonlinear phenomenon as close as possible. These linearizations are of finite dimensions and exceed the phase space dimension of the nonlinear system by one at most. In particular, this shows, that linear immersions of nonlinear system are not necessarily high-dimensional. We introduce simple, low-dimensional embeddings of nonlinear dynamical systems into locally topologically equivalent linear systems, see Fig.~\ref{fig:overview}. More specifically, we investigate nonlinear systems with continuous \mycorr{Koopman} spectrum, limit cycles, and  coexisting attractors. The systems we consider have either eluded linearizations~\cite{lusch2018deep}, have only been partially linearized \cite{korda2018linear,korda2020optimal},  or been proven to be non-linearizable in the sense of topological conjugacy of the same dimension~\cite{cenedese2022data}.

\section{Linear Immersions of Nonlinear Dynamical Systems}
\label{sec:theory}

%In this section, we outline the set-up of our analysis (cf. Section~\ref{sec:setup}) and define the objects of interest. We then specify our learning algorithm and elucidate our approach to low-dimensional linear immersions.

%\subsection{Set-Up}
%\label{sec:setup}
We consider the general, autonomous dynamical system  
\begin{equation}
\label{eq:NL_system}
    \dot{\mathbf{x}}= \mathbf{f}(\mathbf{x}),\qquad \mathbf{x}\in \mathcal{X}\subseteq{\mathbb{R}}^N,
\end{equation}
where $\mathcal{X} \subseteq{\mathbb{R}}^N$ denotes the domain of definition and $\mathbf{f}$ is assumed to be sufficiently smooth to guarantee existence and uniqueness of solutions.  Let $\bm{F}^t(\mathbf{x}_0)$ denote the flow map of system \eqref{eq:NL_system} with the initial condition $\bm{x}_0$.  A \textit{linear immersion} in the sense of \cite{mauroy2020koopman,liu2023properties} is an injective transformation \mycorr{$\bm{\Phi}:\mathcal{X} \to \mathcal{Y}$, for an open set $\mathcal{Y}\subseteq \mathbb{R}^M$,} that maps trajectories of the nonlinear system~\eqref{eq:NL_system} into trajectories of a linear system
    \begin{equation}
\label{eq:lin_system}
    \dot{\mathbf{y}}= \mathbf{A}\mathbf{y}, \qquad \mathbf{y}\in\mathcal{Y}\subseteq{\mathbb{R}}^M,
\end{equation}
with $M\geq N$, such that 
\begin{equation}\label{immserion}
    \bm{\Phi}(\bm{F}^t(\bm{x}_0)) = e^{\bm{A}t}\bm{y}_0,
\end{equation}
holds for all initial conditions  $\bm{x}_0\in \mathcal{X}$, $\bm{y}_0\in\mathcal{Y}$ and all times $t\geq 0$.

The transformation $\bm{\Phi}$ together with its inverse $\bm{\Phi}^{-1}$ (restricted to $\mathcal{Y}$) thus allows us to evaluate trajectories of the nonlinear system \eqref{eq:NL_system} by mapping initial conditions under $\bm{\Phi}$, evolution of the linear dynamics, and transforming back:
\begin{equation}
    \label{eq:x_lin}
    \begin{split}
            \mathbf{F}^{t}(\mathbf{x}_0)&=\bm{\Phi}^{-1}(e^{\mathbf{A}t}\mathbf{y}_0)\\
           & =\bm{\Phi}^{-1}(e^{\mathbf{A}t}\bm{\Phi}(\mathbf{x}_0)),
    \end{split}
\end{equation}
as illustrated in the diagram~\ref{fig:commu}.
\tikzcdset{arrow style=tikz, diagrams={>=stealth'}}
\begin{figure}[h!]
 	\begin{center}
%			\[ \psset{arrows=->, arrowinset=0.25, linewidth=0.6pt, nodesep=3pt, labelsep=2pt, rowsep=1.1cm, colsep = 1.1cm, shortput =tablr}
% \everypsbox{\scriptstyle}
% \begin{psmatrix}
% \mathbf{y}_0 \in \mathcal{Y} & e^{\mathbf{A}t} \mathbf{y}_0\\%
% \mathbf{x}_0 \in\mathcal{X}   & \mathbf{F}^t(\mathbf{x}_0)
 %%%
% \ncline{1,1}{1,2}^{e^{\mathbf{A}t} } 
% \ncline{2,1}{1,1} <{\bm{\Phi} }
% \ncline{1,2}{2,2} > {\bm{\Phi}^{-1}}
% \ncline{2,1}{2,2}^{\mathbf{F}^t}
% \end{psmatrix}
% \]
\[
  \begin{tikzcd}
     \mathbf{y}_0 \in \mathcal{Y} \arrow{r}{e^{\mathbf{A}t}}   & e^{\mathbf{A}t} \mathbf{y}_0 \arrow{d}{\bm{\Phi}^{-1}} \\[0.5cm]
     \mathbf{x}_0 \in\mathcal{X}  \arrow{r}{\mathbf{F}^t} \arrow{u}{\bm{\Phi}} & \mathbf{F}^t(\mathbf{x}_0)
  \end{tikzcd}
\]
 	\end{center}
 	\caption{The mappings $\bm{\Phi}$ and $\bm{\Phi}^{-1}$ allow to generate trajectories of the nonlinear system~\eqref{eq:NL_system} using linear dynamics only.    }
 	\label{fig:commu}
\end{figure}

We stress that the dimension of the linear system $M$  does not have to be the same as the dimension of the nonlinear system $N$. While the linear system in Koopman theory is generally infinite dimensional ($M=\infty$), more classical linearization results  \cite{guckenheimer2013nonlinear} establish the existence of local linearizing transformations between spaces of equal dimensions (i.e. $M=N$). 
%Thus, it is generally anticipated that the linear system~\eqref{eq:lin_system} has at least the dimension of the nonlinear counterpart~\eqref{eq:NL_system}, i.e. $M \geq N$. 

Although the dimension $M$ of the ambient space of the linear dynamics can be larger than the dimension of the nonlinear dynamics in our approach, the dimension of $\mathcal{Y}$ has to be $N$ in order to guarantee dynamical conjugacy of the nonlinear system~\eqref{eq:NL_system} and the linear system~\eqref{eq:lin_system} almost everywhere.

%For the dynamical system~\eqref{eq:NL_system}, we seek a linearization of the forma transformation which transforms trajectories of the nonlinear system~\eqref{eq:NL_system} into to solutions of the linear system~\eqref{eq:lin_system} and the reverse transformation $\Phi^{-1}$ mapping trajectories of the linear system~\eqref{eq:lin_system} to solutions the nonlinear system~\eqref{eq:NL_system}.This corresponds to the commutative diagram is sketched in FIG 1 OR SOMETHING NEW.  This allows to obtain solutions of the nonlinear system~\eqref{eq:NL_system} utilizing linear dynamics as follows 

  %We note that any time-dependent dynamical system $ \mathbf{x}= \mathbf{f}(\mathbf{x},t)$ can be made autonomous by adding time as a dynamic variable with trivial dynamics $\dot{t}=1$. \todo[inline, color = green]{Wenn du meinst, dass wir eine wasserdichte definition brauchen, kannst du den text noch einbauen.}

Established approaches leveraging Koopman theory, such as DMD~\cite{schmid2010dynamic} and its extensions~\cite{williams2016extending} attempt to construct the linearizing transformation ($\bm{\Phi}$ in our notation) and the linear dynamics~\eqref{eq:lin_system} in a one-step procedure. This ambitious approach, however, often results in complex and uninterpretable outcomes even for simple nonlinear systems. 

Here, we follow a different approach by choosing  the linear dynamics~\eqref{eq:lin_system} from an appropriate set of reference systems. These linear dynamics  are either topologically equivalent to the nonlinear dynamics~\eqref{eq:NL_system} or mimic them closely. This simplification allows us to construct low dimensional and easy-to-interpret linear representations of nonlinear dynamics. Moreover, this equivalence implies that there exists a mapping $\bm{\Phi}$ which maps trajectories of the nonlinear system~\eqref{eq:NL_system} into trajectories of the linear system~\eqref{eq:lin_system} and the reverse transformation $\bm{\Phi}^{-1}$ which maps trajectories of the linear system~\eqref{eq:lin_system} into trajectories of the nonlinear system~\eqref{eq:NL_system}. This justification allows us to compute these mappings and, in principle, various methods could be applied to calculate or approximate them. In the following, we either rely on explicit calculations or utilize universal function approximators, i.e. neural networks. 

This approach draws from an extensive literature and decade long research in nonlinear dynamics. Indeed, the principal phase space structure  of nonlinear systems can often be obtained by, e.g., launching trajectories~\cite{nayfeh2008applied}, cell-mapping methods~\cite{hsu2013cell}, spectral analysis~\cite{kaplan2012understanding} or geometric arguments~\cite{guckenheimer2013nonlinear}. We propose to build on such preliminary analyses and select a linear system~\eqref{eq:lin_system} based on these insights. We do not seek to answer whether such an envisioned linearization is achievable in the most general setting, nor do we seek to establish a routine or algorithm to construct such a linearization in general - this exploratory approach might be deferred to a forthcoming work. Rather, by considering explicit examples of nonlinear systems, we exemplify that our approach yields linear representations of nonlinear dynamics featuring continuous \mycorr{Koopman} spectra, limit cycles, and coexisting steady states. 

\subsection{Learning linear immersions}
\label{sec:learning}
Once the conjugated linear dynamics~\eqref{eq:lin_system} have been selected, the transformation $\bm{\Phi}$ and $\bm{\Phi}^{-1}$ can be approximated by deep neural networks. In this work, we use feedforward neural networks with the Tansig activation function, i.e.,  $\mbox{tansig}(x)=2/(1+e^{-2x})-1$, see~\cite{schmidhuber2015deep}. \mycorr{In all our examples, the networks consist of three layers with twenty hidden units each, except for the Van-der-Pol oscillator discussed in Section~\ref{sec:LC}, where we increase the number of hidden units to forty.}

To obtain the parameters of the neural networks, we generate trajectories of the nonlinear system~\eqref{eq:NL_system} and the linear equivalent~\eqref{eq:lin_system} utilizing~\mycorr{a numerical time integration scheme based on the Dormand-Prince method~\cite{dormand1980family}, as implemented in MATLAB's routine \texttt{ode45}.} Collecting the trajectory data, we seek to minimize  
\begin{equation}
\label{eq:minimizations}
 |\mathbf{y}-\mathbf{\Phi}(\mathbf{x})|^2\to\min,\qquad |\mathbf{x}-\mathbf{\Phi}^{-1}(\mathbf{y})|^2\to\min
 \end{equation}
 for all $\mathbf{y}\in \mathcal{Y}$ and $\mathbf{x}\in \mathcal{X}$. We seek approximate solutions to both minimization~\eqref{eq:minimizations}. To this end, we select the Levenberg--Marquardt algorithm, an interpolation between the Gauss--Newton algorithm and gradient descent~\cite{kelley1999iterative}, to optimize the network weights of the feedforward neural networks.~\mycorr{In this work we employ MATLAB's routine \texttt{train} for this optimization.} While the inversion of the Jacobian in the Levenberg--Marquardt algorithm increases the computational costs, we found that the residual errors were considerably smaller compared to a pure gradient descent approach.

 \mycorr{Approximating the mappings $\boldsymbol{\Phi}$ and $\boldsymbol{\Phi}^{-1}$ with deep neural networks introduces errors and hence trajectories on the nonlinear system~\eqref{eq:NL_system} and the transformed trajectories of the linear system~\eqref{eq:lin_system} are not exactly equal. In Appendix~\ref{app:error}, we obtain an upper bound on this difference. The derived upper bound depends linearly on the residuals of the minimization~\eqref{eq:minimizations}, on the growths of the linear system~\eqref{eq:lin_system} and a Lipschitz constant of the neural network approximating the mapping $\boldsymbol{\Phi}^{-1}$. For linear systems with bounded trajectories, this bound remains finite for all times. }

We emphasize that the key property of our method - besides the specific architecture of the neural network - lies in the choice of immersed linear dynamics, i.e., the matrix $\bm{A}$ in \eqref{eq:lin_system}. We consider specific reference systems depending on the nonlinear phenomena we want to embed. As mentioned before, these include systems with continuous \mycorr{Koopman} spectra, limit cycles, and systems with coexisting steady states. 

In addition to the data-driven computations, we provide an overview of explicit embeddings for each of these nonlinear phenomena in Appendix~\ref{app:examples}.

\section{Continuous Koopman spectra: The Nonlinear Pendulum and the Bistable Duffing Oscillator}
\label{sec:continuous}

\mycorr{
In this section, we discuss global embeddings of nonlinear systems with nested families of periodic orbits. These systems are challenging for Koopman embeddings as the Koopman operator is expected to have a continuous operator spectrum \cite{mezic2020spectrum}. More specifically, the corresponding Koopman eigenfunctions are no longer eigenfunctions as in the case of isolated eigenvalues but have to be interpreted in the sense of distributions \cite{hislop2012introduction}.
}

We first present an explicit global, linear representation of class of nonlinear systems featuring a one parameter family of periodic orbits in Appendix \ref{app:harmonic}. Therein, the linear harmonic oscillator is selected as linear system~\eqref{eq:lin_system}. After a transformation into the polar coordinates we obtain a nonlinear system with a nested family of periodic orbits with continuously varying frequency. 

After this explicit construction for specific nonlinear systems, we construct a data-driven linear representation for the two popular benchmark systems: i) the nonlinear pendulum and ii) the \mycorr{bistable}, conservative Duffing oscillator.

\mycorr{For the nonlinear pendulum specifically, Lusch \textit{et al.}~\cite{lusch2018deep} argue that the energy dependent frequency (i.e. backbone curve) and the appearance of higher harmonics give rise to a continuous spectrum of the Koopman operator and claim that a low dimensional linearization is therefore impossible.} To alleviate this problem, an auxiliary network to parameterize the continuous \mycorr{Koopman} spectrum of the nonlinear pendulum is proposed. This approach, however, renders the dynamics in Koopman coordinates to be in in the form of $\bm{y}_{k+1}=\bm{K}(\bm{y}_k)\bm{y}_k$. This system is of course nonlinear. Indeed, the authors point out that their Koopman coordinates reassemble action angle variables \cite{arnol2013mathematical}. The dynamics of the nonlinear pendulum in action-angle coordinates are only linear in the center directions (conserved quantities), while the dependence of the phase on those conserved quantities is in general nonlinear. Thus, the global linearization of the nonlinear pendulum remains incomplete.

\subsection{The Nonlinear Pendulum}
\label{sec:pendulum}

\begin{figure*}

 	\begin{center}
			\includegraphics[width=\textwidth]{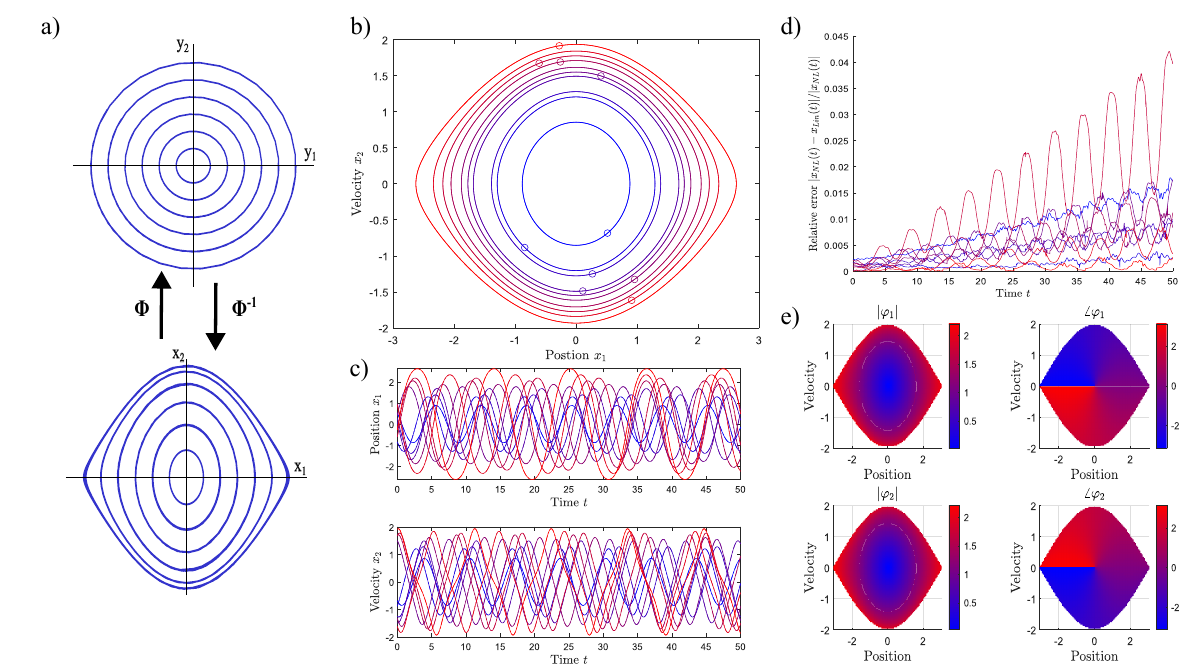}
 	\end{center}
 	\caption{Learning the linear conjugacy of the nonlinear pendulum~\eqref{eq:pendulum}. \textbf{a)} The embedding $\bm{\Phi}$ and its restricted inverse $\bm{\Phi}^{-1}$ map orbits of the nonlinear pendulum into orbits of the linear oscillator and vice versa.  \textbf{b)} Phase space of the transformed linear oscillator~\eqref{eq:lin_pendulum}. Circle markers indicate the randomly selected initial conditions.  \textbf{c)} Raw time series of the position and the velocities of the transformed linear oscillator~\eqref{eq:lin_pendulum}. \textbf{d)} Relative error between trajectories of the transformed linear oscillator~\eqref{eq:lin_pendulum} and the nonlinear pendulum~\eqref{eq:pendulum}.  \textbf{e)} Koopman eigenfunctions of the nonlinear pendulum~\eqref{eq:pendulum}.  }
 	\label{fig:Pendulum}
\end{figure*}

The second-order equation of motion for the nonlinear pendulum is given by
 \begin{equation}
 \label{eq:pendulum}
     \ddot{x}+\sin(x)=0, \quad \iff  \quad 
     \begin{cases}
     \dot{x}_1=x_2,\\
     \dot{x}_2=-\sin(x_1).\\
     \end{cases}
 \end{equation}
The nonlinear pendulum~\eqref{eq:pendulum} is Hamiltonian with the energy $E=x_2^2/2-\cos(x_1)$ and can be integrated explicitly in terms of elliptic functions. For energies $E<1$  the phase space of the pendulum consists of a one parameter family of nested periodic orbits. This geometry is topologically equivalent to the closed trajectories of the linear oscillator
\begin{equation}
\label{eq:lin_osci}
     \ddot{y}+y=0, \quad \iff  \quad 
     \begin{cases}
     \dot{y}_1=y_2,\\
     \dot{y}_2=-y_1.\\
     \end{cases}
\end{equation}
%System \eqref{eq:pendulum} has been studied extensively in the context of Koopman learning, see \cite{lusch2018deep}. It is known that higher harmonics may cause problems in the learning of the Koopman embedding \cite{lusch2018deep}.\\
%Action-angle variables are canonical variables $(\bm{J,\bm{w}})$ for a Hamiltonian system which allows to write the dynamics as\begin{equation}\label{actionangle}    \begin{cases}        \dot{\bm{J}} = 0,\\        \dot{\bm{w}} = \bm{\omega}(\bm{J}).    \end{cases}\end{equation}
In the following, it will be more convenient to work in polar coordinates. Let $(r_{p},\theta_p)$ denote polar coordinates for the phase space of the nonlinear pendulum and let  $(r_{l},\theta_l)$ denote the corresponding polar coordinates for the linear system~\eqref{eq:lin_osci}. Then, the immersion~\eqref{immserion} is written as $(r_l,\theta_l)=\bm{\Phi}(r_{p},\theta_p)$ with the restricted inverse $(r_p,\theta_p)=\bm{\Phi}^{-1}(r_{l},\theta_l)$. These transformations map the phase space of the nonlinear pendulum~\eqref{eq:pendulum} to the phase space of the linear oscillator~\eqref{eq:lin_osci} and vice versa as illustrated in Fig.~\ref{fig:Pendulum}a.

Since each orbit of the pendulum has a different period $T_p(r_{p},\theta_p)$, we do not expect that the family of closed trajectory of the nonlinear pendulum~\eqref{eq:pendulum} can be immersed into the linear oscillator~\eqref{eq:lin_osci} for the fixed linear frequency $T_{\rm Lin}=2\pi$.~\mycorr{Therefore, we also learn the mapping $T_p(r_{p},\theta_p) = T_p(\bm{\Phi}^{-1}(r_{l},\theta_l))$ that relates the phase space coordinates of the linear oscillator~\eqref{eq:lin_osci} to the period of the nonlinear system. To this end, we utilize the same trajectory data and network construction (cf. Section~\ref{sec:learning}) that we employ to learn the mappings $\bm{\Phi}$ and $\bm{\Phi}^{-1}$.} This approach is reminiscent of the frequency modulation in the Poinc\'{e}--Lindstedt series \cite{verhulst2006nonlinear}. 

Consequently, based on the three networks $\bm{\Phi}(r_{l},\phi_l)$, $\bm{\Phi}^{-1}(r_{l},\phi_l)$ and $T_p(r_{l},\phi_l)$, we construct a linearization of the nonlinear pendulum 
\begin{equation}
\label{eq:lin_pendulum}
     \begin{split}
         (r_{\rm Lin}&(t), \theta_{\rm Lin}(t))\\
         &=\bm{\Phi}^{-1}\left(r_l(t),\frac{(\theta_l(t)-\theta_l^0)T_{Lin}}{T_p(r_l,\theta_l)}+\theta_l^0\right),
     \end{split}
\end{equation}
where $\theta_l^0$ is the initial angle.

Fig.~\ref{fig:Pendulum}b-d compares trajectories of the nonlinear pendulum~\eqref{eq:pendulum} with trajectories of the transformed  linear oscillator~\eqref{eq:lin_pendulum} for randomly selected initial conditions. Overall, excellent agreement is observed. The relative error between the transformed linear system~\eqref{eq:lin_pendulum} and the nonlinear pendulum remains below 5 \% (cf. see Fig.~\ref{fig:Pendulum}d).~\mycorr{Moreover, visualizations of the approximated embedding $\bm{\Phi}$ and its restricted inverse $\bm{\Phi}^{-1}$ are included in the Appendix~\ref{app:pendulum_phi_visu}.}

\mycorr{The errors shown in Fig.~\ref{fig:Pendulum}d grow in time. This growth can be explained as follows. First, we utilize the transformation $\bm{\Phi}$ to map initial conditions of the nonlinear pendulum~$\mathbf{x}_0$ to initial conditions of the linear oscillator $\mathbf{y}_0=\bm{\Phi}(\mathbf{x}_0)$. Since we approximate the mapping $\bm{\Phi}$ with a deep neural network (cf. Section~\ref{sec:learning}), we in fact obtain the point $\tilde{\mathbf{y}}_0$ when applying the obtained approximation of $\bm{\Phi}$ to the initial condition $\mathbf{x}_0$. The point $\tilde{\mathbf{y}}_0$ is in the vicinity of the true initial condition $\mathbf{y}_0$ and the distance between these points (i.e. $|\tilde{\mathbf{y}}_0-\mathbf{y}_0|$) decreases with  decreasing residual of the minimization~\eqref{eq:minimizations} (cf. the error bound~\eqref{eq:error_bound} in Appendix~\ref{app:error}). This observation explains the initial error at $t=0$ in Fig.~\ref{fig:Pendulum}d. Initializing the linear oscillator~\eqref{eq:lin_osci} at  $\tilde{\mathbf{y}}_0$  generates a closed orbit $\mathbf{y}(t)$. Mapping this orbit back utilizing an approximation of the mapping $\boldsymbol{\Phi}^{-1}$ yields a closed orbit with the frequency $T_p(\bm{\Phi}^{-1}(r_{l},\theta_l))$. Although this orbit is in the vicinity of the true orbit of the nonlinear pendulum $\mathbf{x}(t)$, its shape and frequency will slightly differ from the true orbit of the nonlinear pendulum. Thus, trajectories of the transformed linear oscillator~\eqref{eq:lin_pendulum} and the nonlinear pendulum  are initially close, but slowly desynchronize due to the slight difference in their frequencies induced by approximating the mapping $T_p(\bm{\Phi}^{-1}(r_{l},\theta_l))$.  Thus, the error grows as seen in Fig.~\ref{fig:Pendulum}d. Similar observations are also made for solutions obtained from classical method of averaging which can guarantee a small error on long, but finite time scales~\cite{verhulst2006nonlinear}. Decreasing the residual of the minimization~\eqref{eq:minimizations} will generally decrease the error shown in Fig.~\ref{fig:Pendulum}d.}

\mycorr{We note that all orbits shown in Fig.~\ref{fig:Pendulum} are bounded and hence the error remains bounded. Indeed, the boundedness of the approximation error can also be deduced from the upper bound~\eqref{eq:error_bound} and by noting that the flow map of linear oscillator~\eqref{eq:lin_osci} in bounded.}

We stress that normalizing the angular coordinate $\theta_l(t)$ by the learned period $T_p(r_l,\theta_l)$ in the transformation~\eqref{eq:lin_pendulum} is crucial in order to obtain a global immersion. Through this normalization the fixed period of the harmonic oscillator is mapped to the correct energy-dependent period of the nonlinear pendulum. Thus, each orbit shown in Fig.~\ref{fig:Pendulum}b is traversed within the period of the nonlinear pendulum and the correct raw time series are reconstructed, see Fig.~\ref{fig:Pendulum}c. Without this normalization the mapping $\bm{\Phi}^{-1}$ would correctly reconstruct the nonlinear phase space geometry shown in Fig.~\ref{fig:Pendulum}b, but would fail to capture the raw time series of the coordinates $x_1$ and $x_2$ (which would have the period of the linear pendulum $T_{\rm Lin}$). To illustrate the necessity to include the nonlinear period $T_p(r_l,\theta_l)$ in the transformation~\eqref{eq:lin_pendulum}, we show trajectories of the transformed linear oscillator without frequency adjustment in Fig.~\ref{fig:Pendulum_noTcorr} in Appendix~\ref{app:Tcorr}. While the phase space geometry is correctly reconstructed (cf. Fig.~\ref{fig:Pendulum_noTcorr}a), the time series of the linear oscillator without frequency adjustment~\eqref{eq:lin_pendulum_noTcorr} deviate quickly from the trajectories of the nonlinear pendulum (cf. Fig.~\ref{fig:Pendulum_noTcorr}c). We note that polar coordinates altogether facilitate the implementation of this correction mechanism in the transformation~\eqref{eq:lin_pendulum}. Thereby, we avoid the nonlinear dynamics in the Koopman coordinates proposed in~\cite{lusch2018deep}. 

Transforming the linear oscillator~\eqref{eq:lin_osci} into a diagonal system yields the coordinates ($\varphi_1,\varphi_2$) which define Koopman eigenmodes of the nonlinear pendulum~\eqref{eq:pendulum} (cf. Appendix~\ref{app:Koopman} for a brief review of the relationship between linear immersions and the Koopman operator). Their magnitude and phase are shown in Fig.~\ref{fig:Pendulum}e. Since the choice of the diagonalization is only unique up to a scalar multiple of the eigenvectors, the moduli $|\varphi_1|$ and $|\varphi_2|$ are only determined up to an arbitrary constant. Similarly, the angles $\angle \varphi_1$ and $ \angle \varphi_2$ can be rotated by an arbitrary angle. The main features of the Koopman eigenmodes, however, remain invariant under these symmetries. Their respective moduli grow with the distance to the origin. Moreover, the angle of the first Koopman eigenmode  $\varphi_1$ (corresponding to the eigenvalue with a positive imaginary part) increases in the clockwise direction, whereas the angle of the second Koopman eigenmode $\varphi_2$ increases in the counter-clockwise direction.

The linear oscillator~\eqref{eq:lin_osci} with frequency equal to one can be interpreted as a particular two-dimensional subspace of the infinite-dimensional Koopman operator. Indeed, any linear oscillator of the form
\begin{equation}
\label{eq:lin_osci_ws}
     \ddot{y}+\omega^2 y=0, \quad \iff  \quad 
     \begin{cases}
     \dot{y}_1=y_2,\\
     \dot{y}_2=-\omega^2y_1,\\
     \end{cases}
\end{equation}
for any frequency $\omega \in \mathbb{R} \! \setminus \!\{0\}$ is equivalent to the nonlinear pendulum via the transformation
 \begin{equation}
\label{eq:lin_pendulum_trans}
     \begin{split}
         (r_{\rm Lin} & (t), \theta_{\rm Lin}(t)) \\
         & =\bm{\Phi}^{-1}\left(r_l(t),\frac{ (\theta_l(t)-\theta_l^0) 2\pi}{\omega  T_p(r_l,\theta_l)}+\theta_l^0\right).
     \end{split}
\end{equation}
Thus, the spectrum of the Koopman operator includes the whole imaginary line as discussed in~\cite{mezic2020spectrum}. A direct approximation of infinite-dimensional Koopman embeddings with neural network is expected to be not only computationally expensive but also lacks a general, theoretical basis. The universal function approximation theorem - the theoretical underpinning of most neural network approaches - is only formulated for mappings between finite dimensional spaces \cite{pinkus1999approximation}. Hence, a direct approximation of the full Koopman operator is generally out of reach and finite-dimensional, discretized version are used instead. We avoid these complications by restricting the transformations $\bm{\Phi}$ and $\bm{\Phi}^{-1}$ to a two-dimensional subspace. Nevertheless, we can deduce the full, continuous spectrum of the Koopman operator.

\begin{figure}[h!]
 	\begin{center}
			\includegraphics[width=0.5\textwidth]{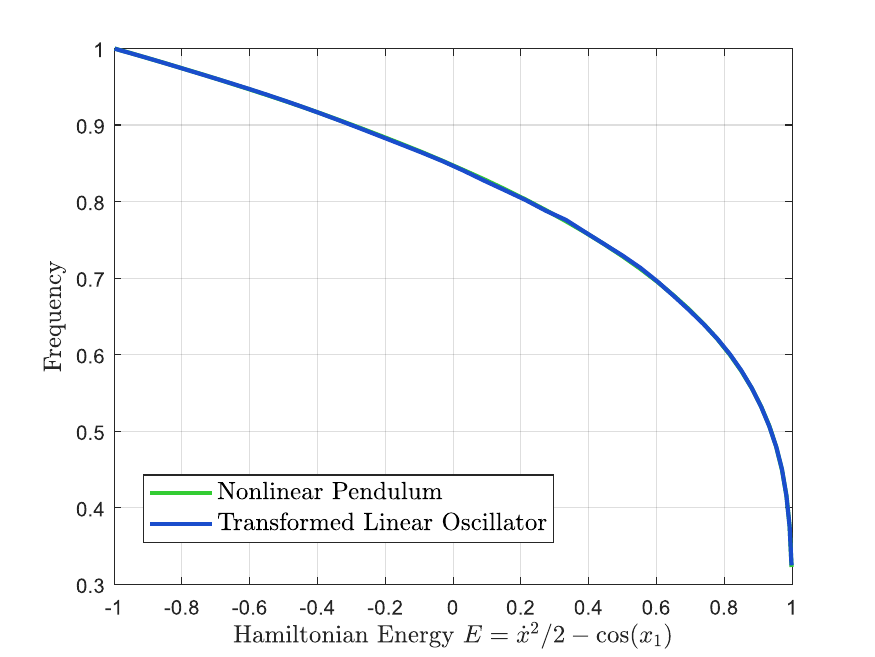}
 	\end{center}
 	\caption{Nonlinear frequency-energy relationship of the nonlinear pendulum~\eqref{eq:pendulum} and the transformed linear oscillator~\eqref{eq:lin_pendulum}. }
 	\label{fig:pendulum_bakcbone}
\end{figure}

In addition, the nonlinear frequency as a function of the energy can be extracted from the transformed linear oscillator~\eqref{eq:lin_pendulum}. \mycorr{This frequency-energy dependency is often utilized to illustrate nonlinear normal modes, a popular tool to analyze vibrations of nonlinear structures~\cite{kerschen2009nonlinear}. For more details on nonlinear normal modes and their use in nonlinear modal analysis, we refer to~\cite{kerschen2014modal}.} We depict this relationship, extracted from the transformed linear oscillator~\eqref{eq:lin_pendulum} along with the curve extracted from the nonlinear pendulum~\eqref{eq:pendulum} in Fig.~\ref{fig:pendulum_bakcbone}. Both curves are practically indistinguishable.

In this section, we have limited our discussion to energy levels below one, for which the phases space consists of a one-parameter family of periodic orbits. For higher energy levels, two other one-parameter families of periodic orbits emerge in the phase space $[0,2\pi] \times \mathbb{R}$. In the next section, we consider the Duffing oscillator as an example with multiple coexisting one-parameter families of periodic orbits.

\subsection{\mycorr{Bistable} Duffing Oscillator}
\label{Duffing}

\begin{figure*}
 	\begin{center}
			\includegraphics[width=\textwidth]{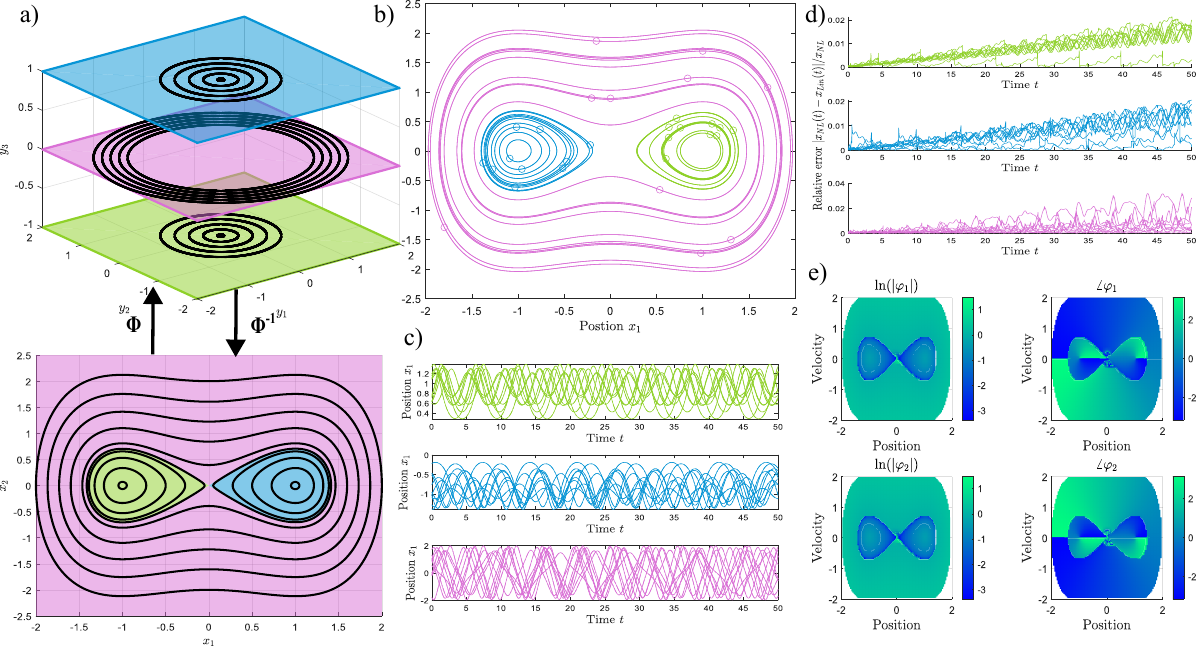}
 	\end{center}
 	\caption{Linearization of the \mycorr{bistable,} conservative Duffing oscillator~\eqref{eq:Duff_cons}. \textbf{a)} Introducing the additional observable $x_3$ with trivial dynamics ($\dot{x}_3 = 0$) allows us to map the different orbit of the conservative Duffing~\eqref{eq:Duff_cons} to the extended linear oscillator~\eqref{eq:ext_lin_osci}.  \textbf{b)} Phase space of the transformed, extended linear oscillator~\eqref{eq:lin_Duff}. The circle markers indicate randomly selected initial conditions. \textbf{c)} Raw time series of the position coordinate $x_1$ for the three distinct phase space domains. \textbf{d)} Relative error between trajectories of the transformed, extended linear oscillator~\eqref{eq:lin_Duff} and the \mycorr{bistable,} conservative Duffing oscillator~\eqref{eq:Duff_cons}. \textbf{e)} Two Koopman eigenfunctions of the \mycorr{bistable} Duffing oscillator~\eqref{eq:Duff_cons}. }
 	\label{fig:Duffing}
 \end{figure*}
 
The \mycorr{bistable,} conservative Duffing oscillator is a Hamiltonian system with energy $E(x_1,x_2) = x_2^2/2-x_1^2/2+x_1^4/4$ and dynamics
\begin{equation}
 \label{eq:Duff_cons}
     \ddot{x}-x+x^3=0, \quad \iff  \quad 
     \begin{cases}
     \dot{x}_1=x_2,\\
     \dot{x}_2=x_1-x_1^3.
     \end{cases}
 \end{equation}
The phase space of the \mycorr{bistable} Duffing oscillator~\eqref{eq:Duff_cons} consists of three nested families of periodic orbits as shown in the bottom on Fig.~\ref{fig:Duffing}a. Two families of periodic orbits emerge from the fixed points $x_1=\pm\sqrt{2}$. These orbits correspond to intrawell oscillations. Furthermore, another family of periodic orbits exists for energies $E>0$. Each of these three families of periodic orbits can be transformed into a linear oscillator as carried out in the preceding Section~\ref{sec:pendulum}. To account for the three families of periodic orbits, we introduce a third observable $x_3$ that encodes the phase space region:
\begin{equation}\label{x3}
\begin{cases}
     x_3 = 1 &\text{ for } E<0, x_1>0,\\
     x_3 = -1 &\text{ for } E<0, x_1<0,\\
     x_3 = 0 &\text{ for } E>0.
\end{cases}
\end{equation}
Since the observable \eqref{x3} is constant along trajectories, we can transform the \mycorr{bistable,} conservative Duffing equation into the extended linear oscillator of the form
\begin{equation}
 \label{eq:ext_lin_osci}
 \begin{split}
      \dot{y}_1&=y_2,\\
     \dot{y}_2&=-y_1,\\
     \dot{y}_3&=0.\\
     \end{split}
 \end{equation}
Similar to the nonlinear pendulum in Section~\ref{sec:pendulum}, we consider the dynamics \eqref{eq:Duff_cons} and \eqref{eq:ext_lin_osci} in polar coordinates. The polar coordinates of the \mycorr{bistable}  Duffing oscillator~\eqref{eq:Duff_cons} are denoted by $(r_D,\theta_D)$, whereas  $(r_l,\theta_l)$ denote the polar coordinates of the first two coordinates of the extended, linear oscillator~\eqref{eq:ext_lin_osci}. This notation, allows us to approximate the embedding $\bm{\Phi}(r_D,\theta_D,x_3)$, its restricted inverse $\bm{\Phi}^{-1}(r_l,\theta_l,y_3)$ and the period $T_D(r_l,\theta_l,y_3)$ by deep neural networks. Then, the linear dynamics take the form 
  \begin{equation}
\label{eq:lin_Duff}
     \begin{split}
         &(r_{\rm Lin} (t),\,  \theta_{\rm Lin}(t)) \\
         &\quad = \bm{\Phi}^{-1}\left(r_l(t),\frac{(\theta_l(t)-\theta_l^0)T_{\rm Lin}}{T_D(r_l,\theta_l,y_3)}+\theta_l^0,y_3\right), 
     \end{split}
\end{equation}
where $T_{\rm Lin}=2\pi$ denotes the period of the extended linear oscillator~\eqref{eq:ext_lin_osci}. Fig.~\ref{fig:Duffing}b depicts the phase space of the transformed, extended linear oscillator~\eqref{eq:lin_Duff}. These  trajectories are practically indistinguishable from the trajectories of the nonlinear Duffing oscillator~\eqref{eq:Duff_cons}. The time series of the position show excellent agreement as well (cf. Fig.~\ref{fig:Duffing}c) and the relative error of trajectories of the transformed, extended linear oscillator remains below two percent, see Fig.~\ref{fig:Duffing}d.~\mycorr{As discussed in the previous section (cf. Section~\ref{sec:pendulum}) this error grows due to the slow desynchronization between the true orbit of the bistable Duffing oscillator~\eqref{eq:Duff_cons} and the orbit of the transformed, extended linear oscillator~\eqref{eq:lin_Duff}.}

Additionally, we compute the two nontrivial Koopman eigenfunctions of the extended linear oscillator~\eqref{eq:ext_lin_osci} and depict their magnitude and phase in Fig.~\ref{fig:Duffing}e. The three distinct phase regions are clearly discernible. As discussed in the previous section~\ref{sec:pendulum}, these Koopman eigenfunctions only span finite-dimensional invariant subset of the infinite dimensional Koopman operator. Selecting a different oscillation frequency in the extended linear oscillator~\eqref{eq:ext_lin_osci} yields other spectral components of the infinite dimensional Koopman operator.

\begin{figure}[h!]
 	\begin{center}
			\includegraphics[width=0.5\textwidth]{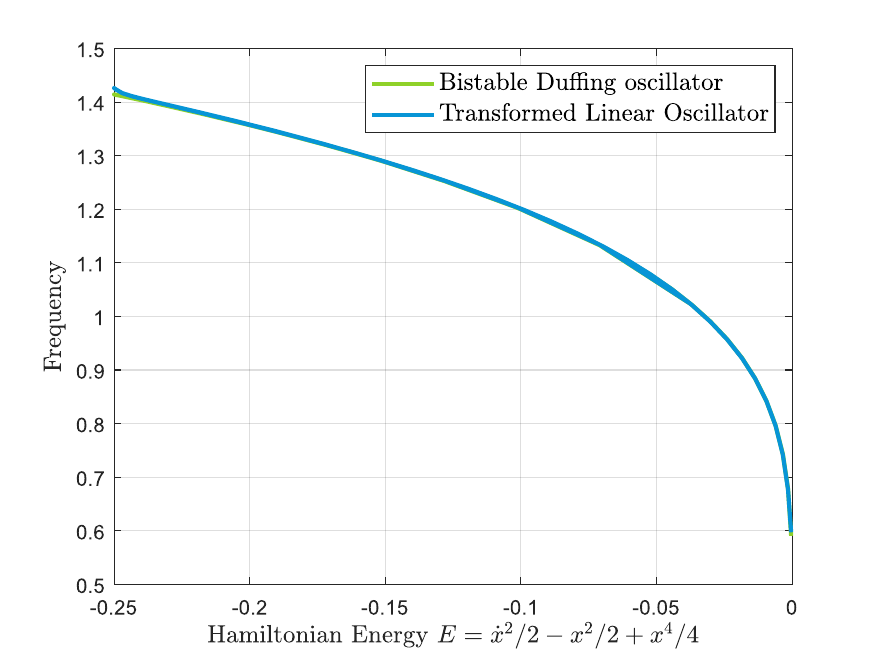}
 	\end{center}
 	\caption{Nonlinear frequency-energy relationship of families of periodic orbits emanating from the center fixed points at $\pm\sqrt{2}$.  }
 	\label{fig:duf_bakcbone}
\end{figure}

Furthermore, the frequency energy relationship for the periodic orbits emanating from the center fixed points $x_1=\pm\sqrt{2}$ are extracted from the extended linear oscillator~\eqref{eq:ext_lin_osci} and the nonlinear Duffing oscillator~\eqref{eq:Duff_cons} are shown in Fig.~\ref{fig:duf_bakcbone}. Both curves are in excellent agreement. 

The crucial extension presented in this section compared to the nonlinear pendulum lies in the introduction of the dynamically trivial observable $x_3$ in equation~\eqref{x3}. This observable separates the three one-parameter families of periodic orbits in the extended phase space (cf. Fig.~\ref{fig:Duffing}a). Each family of periodic orbits is contained in a different level set $\{x_3= \rm const.\}$.

This allows us to construct a system that is topologically equivalent to the extended linear oscillator~\eqref{eq:ext_lin_osci} as follows. First, consider the periodic orbits in any of the three planes $x_3=\{-1,0,1\}$. These families can be smoothly continued until the corresponding family covers the entire plane $\mathbb{R}^2$. In the $x_3$-direction, we can choose an arbitrary continuous continuation of the families of periodic orbits, such that they match with those of the \mycorr{bistable}  Duffing oscillator at the three level sets $x_3=\{-1,0,1\}$. This yields the phase space of an extended Duffing oscillator, which consists of planes $x_3= \rm const.$ foliated by one parameter families of periodic orbits. This phase space geometry is topologically equivalent to the phase space of the extended linear oscillator~\eqref{eq:ext_lin_osci}. This topological equivalence allows us to approximate the corresponding embedding with neural networks. To reduce the computational costs, we restrict the learning of the embedding to the phase space regions corresponding to the solutions of the \mycorr{bistable} Duffing oscillator~\eqref{eq:Duff_cons}. 

\mycorr{To successfully execute this strategy, it is crucial to realize that the dynamics of the bistable Duffing oscillator~\eqref{eq:Duff_cons} consists of three families of periodic orbits. Attempts to embed the dynamics of the bistable Duffing oscillator~\eqref{eq:Duff_cons} to the linear oscillator~\eqref{eq:lin_osci} are bound to fail, since no continuous homomorphism between the phase space of the bistable Duffing oscillator~\eqref{eq:Duff_cons} with three families of periodic obits to the phase space of the linear oscillator~\eqref{eq:lin_osci} with a single family of periodic orbits exists. This further illustrates the necessity to understand the phase space geometry of the nonlinear dynamics first and then subsequently obtain a linearizing transformation as advocated in Section~\ref{sec:theory}.}

The aforementioned technique can be generalized as follows. Assume that we have found linearizing transformations $\{\bm{\Phi}_j(\mathbf{x})\}_{1\leq j\leq n}$ for a certain phase space regions. Then, we introduce an additional coordinate $x^*$ with trivial dynamics $\dot{x}^*=0$ and bundle all the individual transformations into a single embedding $\bm{\Phi}(\mathbf{x},x^*)$. This shows that we can focus on separately linearizing different phase space regions and subsequently patch them together.

Alternative strategies, such as utilizing basis functions with limited support~\cite{williams2015data} or non-smooth transformations~\cite{bevanda2021koopman}, seem less straightforward to approximate with neural networks, which tend to be continuous and supported on the whole domain. Thus, we expect this procedure to be a powerful and versatile approach in obtaining linearizing transformations of complex, nonlinear systems.

%\mycorr{We remark that the two homoclinic orbits of the Duffing oscillator are excluded from our analysis. The measure of this set, however, is zero and hence those will generally not be observed by picking random initial conditions.}  

%If one were to insist to include these two special cases, one could add yet another two level set $y_3=const.$, add the two homoclinic orbits to the training set, and learn an extended embedding. Thereby, one could reproduce any trajectory initialized on the homolcinic orbits. 
\section{Limit Cycles: The Van-der-Pol Oscillator}
\label{sec:LC}

\begin{figure*}[ht!]
 	\begin{center}
			\includegraphics[width=\textwidth]{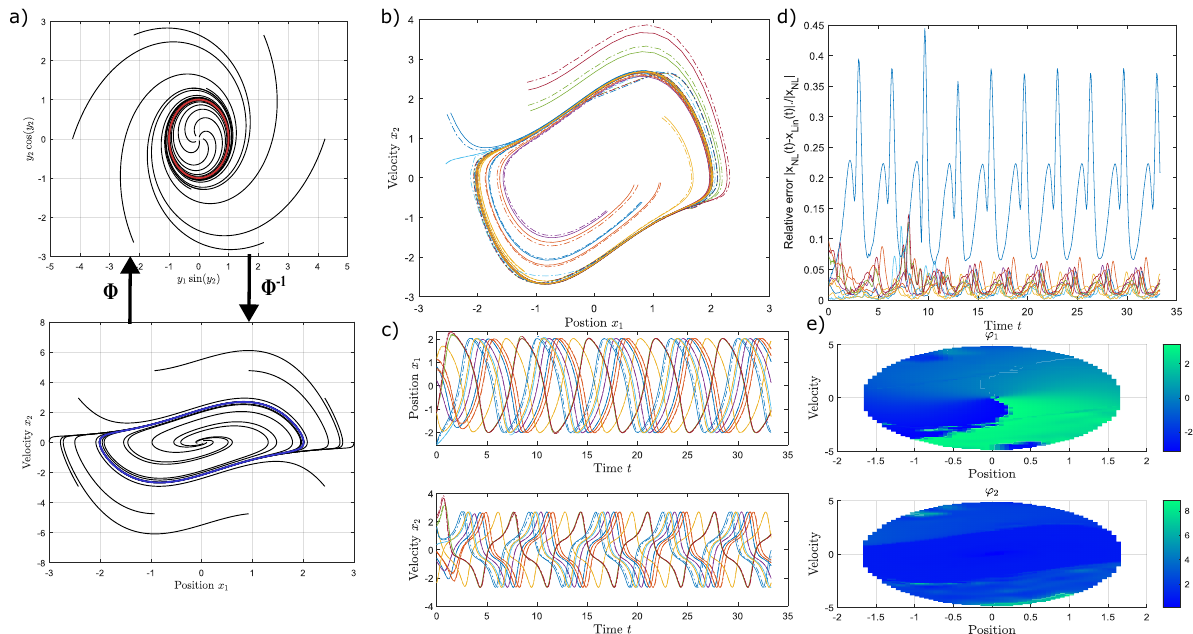}
 	\end{center}
 	\caption{Learning a linear conjugacy for the Van-der-Pol oscillator~\eqref{eq:VdP} \textbf{a)} The embedding $\bm{\Phi}$ and its restricted inverse $\bm{\Phi}^{-1}$ map orbits of the Van-der-Pol oscillator~\eqref{eq:VdP} into orbits of the linear dynamics with limit cycle~\eqref{eq:lin_LC} and vice versa. \mycorr{ \textbf{b)} Phase space of the transformed linear system~\eqref{eq:lin_VdP} (dash-dotted lines) and the Van-der-Pol oscillator~\eqref{eq:VdP} (solid lines).}  \textbf{c)} Raw time series of the position and the velocities of the transformed linear system~\eqref{eq:lin_VdP}. \textbf{d)} Relative error between trajectories of the transformed linear system~\eqref{eq:lin_VdP} and the Van-der-Pol oscillator~\eqref{eq:VdP}.  \textbf{e)} Koopman eigenfunctions of the Van-der-Pol oscillator~\eqref{eq:VdP}. }
 	\label{fig:VdP}
 \end{figure*}

In this section, we learn a global linearization for a planar system with an attracting limit cycle, which constitutes another inherently nonlinear phenomenon.

For a prototypical example showing the conjugacy of a planar nonlinear system with a limit cycle, we refer to Appendix~\ref{app:LC}. Therein, we construct an explicit embedding between a planar nonlinear system with a limit cycle into a three-dimensional linear system. 

As a prototype nonlinear system exhibiting a stable limit cycle we consider the Van-der-Pol oscillator 
\begin{equation}
\begin{split}
    \label{eq:VdP}
     \ddot{x} & +(x^2-1)\dot{x}+x=0, \\
     & \iff  \quad 
     \begin{cases}
     \dot{x}_1=x_2,\\
     \dot{x}_2=-(x_1^2-1)x_2-x_1. 
     \end{cases} 
\end{split}
 \end{equation}
Various attempts to linearize the Van-der-Pol oscillator have been reported in the literature. However, these approaches are either geared to capture the stable limit cycle (e.g., ~\cite{huang2018data,abraham2019active,arbabi2017ergodic}) or the transient dynamics, explicitly excluding the stable limit cycle~\cite{korda2018linear,korda2020optimal}. To the best of our knowledge no global linearization results containing the limit cycle of the Van-der-Pol oscillator have been reported.

We embed the Van-der-Pol oscillator~\eqref{eq:VdP} into the linear system
\begin{equation}
 \label{eq:lin_LC}
     \begin{split}
     \dot{y}_1 & = -D(y_1-1),\\
     \dot{y}_2 & = \Omega,\\
     \end{split}
 \end{equation}
where the coordinate $y_1$ corresponds to the radial coordinate and $y_2$ to the phase. Since the period of the normalized Van-der-Pol limit cycle is given by $T_{\rm VdP}=6.66$, we choose $\Omega=2\pi/T_{\rm VdP}=0.943$ for the linear system~\eqref{eq:lin_LC}. Based on a Floquet analysis we set $D=1.06$, such that the Van-der-Pol oscillator~\eqref{eq:VdP} and the linear system~\eqref{eq:lin_LC} have the same linearized decay rates towards the limit cycle. Approximating the embedding $\bm{\Phi}$ and its restricted inverse $\bm{\Phi}^{-1}$ by deep neural networks, the linear dynamics take the form 
\begin{equation}
\label{eq:lin_VdP}
\begin{split}
(x_1^L(t),& x_2^L(t))\\
&=\bm{\Phi}^{-1}\Big(y_1(t)\sin(y_2(t)),y_1(t)\cos(y_2(t))\Big).
\end{split}
\end{equation}
The trajectories of the transformed linear system~\eqref{eq:lin_VdP} are shown as dash-dotted lines in the Fig.~\ref{fig:VdP}b and Fig.~\ref{fig:VdP}c. While the overall reconstruction is very accurate,  some minor discrepancies remain. For a single trajectory shown in blue in Fig.~\ref{fig:VdP}b and Fig.~\ref{fig:VdP}c, the solution of the Van-der-Pol oscillator~\eqref{eq:VdP} and the transformed linear system~\eqref{eq:lin_VdP} do not synchronize. This asynchronous behavior is also reflected in the relative error (cf., Fig.~\ref{fig:VdP}d). The raw time series in Figure~\ref{fig:VdP}c reveals that the solution of the transformed, linear system~\eqref{eq:lin_VdP} displays a minimal time lag. For all the other initial conditions, however, the relative error remains below 10 percent.~\mycorr{To obtain the results shown in Fig.~\ref{fig:VdP} the number of hidden units is increased to forty (from initially twenty). For twenty hidden units the relative error remains large. Thus, we anticipate that the approximation can be reduced further by employing more powerful neural networks and larger training sets. Additionally, in Appendix~\ref{app:VdP_phi_visu} we show visualizations of the approximated embedding $\bm{\Phi}$ and its restricted inverse $\bm{\Phi}^{-1}$.}

\mycorr{We note that the linear system~\eqref{eq:lin_LC} is an inhomogeneous ordinary differential equation, while the linear system~\eqref{eq:lin_system} is a homogeneous ordinary differential equation. To reformulate the dynamics~\eqref{eq:lin_LC} into a homogeneous system, we may introduce an additional coordinate $y_3$ with trivial dynamics, i.e. $\dot{y}_3=0$. This allows us to rewrite the dynamics~\eqref{eq:lin_LC} as $\dot{y}_1  = -D(y_1-y_3)$ and $\dot{y}_2 = \Omega y_3$. In this setting we restrict our analysis to the invariant subspace $y_3=1$. Thus, the inhomogeneous linear system~\eqref{eq:lin_LC} can be rewritten in the form of the homogeneous linear system~\eqref{eq:lin_system}. }

\mycorr{We select the angular frequency $\Omega$ and the decay rate $D$ of the linear system~\eqref{eq:lin_LC} based on  fundamental characteristics of the Van-der-Pol oscillator~\eqref{eq:VdP}. Generally, the frequency of an limit cycle can be obtained via spectral analysis and decay rates may be estimated by studying the local decay towards the limit cycle. If a different frequency $\Omega$ for the linear system~\eqref{eq:lin_LC} is selected, then the transformed limit cycle of the linear system~\eqref{eq:lin_LC} would slowly desynchronize from the limit cycle of the Van-der-Pol oscillator~\eqref{eq:VdP}, similar to the desynchronization observed for the nonlinear pendulum (cf. Appendix~\ref{app:Tcorr}). For slight variations in the decay rate $D=1$, we obtain residuals for the minimizations~\eqref{eq:minimizations} which are similar to the residuals obtained when selecting $D=1.06$. However, for larger deviations (i.e. $D=0.5$ and $D=2$), the fitting residuals increase significantly and no meaningful transformations~$\bm{\Phi}$ and~$\bm{\Phi}^{-1}$ can be approximated.}

Overall, our construction shows that the Van-der-Pol oscillator~\eqref{eq:VdP} can be transformed into the linear system~\eqref{eq:lin_LC} with a very good accuracy. This prototypical construction can be expanded further to linearize the dynamics about limit cycles of more complex nonlinear systems. To this end, the frequency $\Omega$ of the linear system~\eqref{eq:lin_LC} has to be adapted accordingly. Furthermore, for higher dimensional systems, more decaying coordinates need to be added to the overall dynamics. These decay rates can be chosen, e.g., based on a Floquet analysis of the corresponding limit cycle. Similarly, an extension to quasi-periodic tori can be envisioned. To this end, the linear system~\eqref{eq:lin_LC} needs to be extended by the phase equations $\dot{y}_{k}=\Omega_k$, where $\Omega_k$ denote the frequencies of the quasi-periodic torus. It is anticipated that mirroring the local phase space geometry around the limit cycle simplifies approximating the transformations $\bm{\Phi}$ and $\bm{\Phi}^{-1}$.

\section{Coexisting steady states: The forced-damped Duffing oscillator}
\label{sec:FDDuffing}

\begin{figure*}[ht!]
 	\begin{center}
			\includegraphics[width=\textwidth]{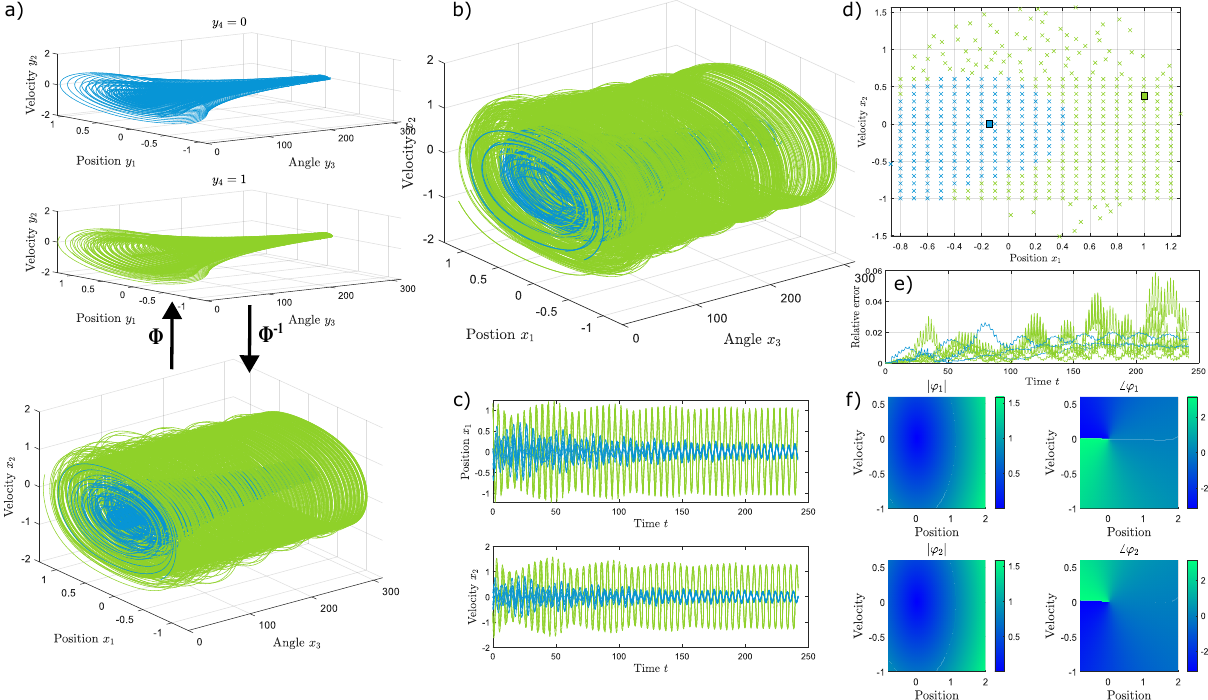}
 	\end{center}
 	\caption{Linearization of the forced-damped Duffing oscillator~\eqref{eq:Duff_forced+damped} with coexisting steady state solutions. 
  \textbf{a)} We learn the embedding between the phase spaces of the forced-damped Duffing oscillator~\eqref{eq:Duff_forced+damped} and the linear system~\eqref{eq:Lin_fd}. 
  \textbf{b)} Testing trajectories of the transformed linear system~\eqref{eq:duff_fd_lin}.  
  \textbf{c)} Time series of the transformed linear system~\eqref{eq:duff_fd_lin}. 
  \textbf{d)} Approximated inflowing invariant phase space region.
  \textbf{e)} Relative error between testing trajectories of the nonlinear Duffing oscillator~\eqref{eq:Duff_forced+damped} and the linear system~\eqref{eq:duff_fd_lin}. \textbf{f)} Koopman eigenfunctions of the forced-damped Duffing oscillator~\eqref{eq:Duff_forced+damped}.}
 	\label{fig:Duff_fd}
 \end{figure*}
 
Systems with coexisting steady states make up another class of nonlinear systems that are difficult to immerse globally. Page and Kerswell~\cite{page2019koopman} indicate that nonlinear systems with multiple invariant solutions cannot be represented globally by a convergent Koopman linearization. The coexistence of multiple limit sets lead Cenedese et al.~\cite{cenedese2022data} to declare the corresponding dynamics as non-linearizable. Moreover, Liu et al.~\cite{liu2023properties} show that linearizing transformation of nonlinear systems with multiple limit sets cannot be one-to-one.

For a prototypical example showing the conjugacy of a planar nonlinear system with multiple fixed-points to a linear system, we refer to Appendix~\ref{app:coexist}. In this model problem, we illustrate how four fixed points (two stable, two unstable) can be represented by the dynamics of a linear system with only \textit{one} globally stable fixed point. 

In the following, we consider the forced-damped Duffing oscillator
\begin{equation}
\label{eq:Duff_forced+damped}
\begin{split}
        \ddot{x} & +c\dot{x}+kx+k_3x^3=f \cos(\Omega t),
    \\[0.33cm]
    &\iff  \quad      \begin{cases}     \dot{x}_1=x_2,\\     \dot{x}_2=-cx_2-kx_1-k_3x_1^3+f \cos(x_3),\\     \dot{x}_3=\Omega ,            \end{cases}
\end{split}
\end{equation}
which constitutes a prototypical example of a nonlinear dynamical system that exhibits multiple limit sets as well as chaotic behavior. Equation~\eqref{eq:Duff_forced+damped} has been used to model forced-damped vibrations of structures~\cite{nayfeh2008nonlinear}, \mycorr{liquid-sloshing}~\cite{bauerlein2021phase}, and hysteresis effects of response curves~\cite{kovacic2011duffing}.

A typical property of the forced-damped Duffing oscillator~\eqref{eq:Duff_forced+damped} is the coexistence of two stable steady state solutions (bi-stability). Indeed, for the parameter values $c=0.02$, $k=k_3=1$, $f=0.1$ and $\Omega=1.3$, trajectories decay towards two distinct stable periodic orbits, as shown the bottom of Fig.~\ref{fig:Duff_fd}a. The nontrivial Floquet multipliers of both stable limit cycles are complex conjugates indicating that the periodic orbits are stable nodes.  

%\todo[inline, color = green]{We need to oint out that we do not approximate a homeomorphsim of the same dimension of the phase space, but a higher dimensional one. This }

For the linear counterpart of the forced-damped Duffing oscillator~\eqref{eq:Duff_forced+damped} we select the damped harmonic oscillator with two additional components:
\begin{equation}
\label{eq:Lin_fd}
    \ddot{y}+c_l\dot{y}+\omega_{l}^2y=0 \quad \iff  \quad 
     \begin{cases}
     \dot{y}_1=y_2,\\
     \dot{y}_2=-c_ly_2-\omega_{l}^2y_1,\\
     \dot{y}_3=\Omega,\\
     \dot{y}_4=0.\\
     \end{cases}
\end{equation}
\mycorr{The coordinate $y_3$ in equation~\eqref{eq:Lin_fd} mimics the period of a limit cycle. Moreover, we utilize the coordinate $y_4$ to account for the two distinct steady states similar to the different periodic regions of the conservative, unforced Duffing oscillator discussed in Section~\ref{Duffing}. More specifically, we restrict the dynamics~\eqref{eq:Lin_fd} to the two level sets $y_4=1$ and $y_4=0$ and associate each of these level sets with the basin of attractions of the stable steady state shown in Fig.~\ref{fig:Duff_fd}a.} We emphasize that we embed the nonlinear dynamics in higher-dimensional space - this allows us to linearize the essentially nonlinear dynamics with coexisting steady states.

Based on  Floquet analysis, we select $c_l=c=0.02$, which yields the same local decay to steady state in the forced-damped Duffing oscillator~\eqref{eq:Duff_forced+damped} and the linear system~\eqref{eq:Lin_fd}. To obtain a value for the linear eigenfrequency $\omega_l$ of the linear system~\eqref{eq:Lin_fd}, we select initial conditions close to the two stable limit cycles of the forced-damped Duffing oscillator~\eqref{eq:Duff_forced+damped}. From the zero-upcrossings of the arising decaying oscillations, we obtain the dominant frequency of $1.3 \mbox{ rad/s}$, and hence, we select  $\omega_{l}=1.3$.

For each level set $y_4=\rm const.$, the phase space geometry of the system~\eqref{eq:Lin_fd} is topologically equivalent to phase space of system~\eqref{eq:Duff_forced+damped} restricted to the basin of attraction of one of the two attractors. Thus, we use the two level sets $y_4=0$ and $y_4=1$ to incorporate the two different steady states. For trajectories with initial conditions decaying to the orbit with lower amplitude (blue in Fig~\ref{fig:Duff_fd}a), we set the coordinate $y_4$ to zero, while we set $y_4=1$ for trajectories decaying to the high amplitude orbit (green in Fig~\ref{fig:Duff_fd}a).

We restrict our analysis to the approximation of an  inflowing (in the sense of \cite{wiggins1994normally}) invariant region shown in Fig.~\ref{fig:Duff_fd}d. To obtain this region, we sample the rectangle $(y_1,y_2)=[-0.8, 1.2]\times [-1,0.6]$, which, although not invariant, contains both attractors. Hence, each trajectory launched inside this rectangle will eventually settle inside the rectangle. To obtain an approximately inflowing invariant region, we track the initial conditions launched at the regular grid points. If those exit the domain and the phase space location is more than $0.1$ away from the closest grid point, we add this phase space location. This yields the inflowing invariant region shown in Fig.~\ref{fig:Duff_fd}d. These state space locations are included as initial conditions in the training data and trajectories were generated for one period. Based on this data, we employ deep neural networks to learn the transformations $\bm{\Phi}$ and $\bm{\Phi}^{-1}$. 

As test data, we randomly select twenty initial conditions $x_1^0$ and $x_2^0$ inside the rectangle shown in Fig.~\ref{fig:Duff_fd}d. Transforming these initial conditions into the phase space of the linear system~\eqref{eq:Lin_fd}, yields the initial conditions $(y_1^0,y_2^0,0,y_4)=\bm{\Phi}(x_1^0,x_2^0,0)$. Subsequently the linear dynamics~\eqref{eq:Lin_fd} yields the trajectories $(y_1(t),y_2(t),y_3(t),y_4(t))$ which are mapped back via
\begin{equation}
\label{eq:duff_fd_lin}
    (x_1^L,x_2^L,x_3^L)=\bm{\Phi}^{-1}(y_1,y_2,y_3,y_4).
\end{equation}
We include the transformed trajectories~\eqref{eq:duff_fd_lin} in Figure~\ref{fig:Duff_fd}. We observe that the transients as well as the steady state behavior is correctly predicted by the transformed linear system~\eqref{eq:duff_fd_lin}. In the phase space, see~Figure~\ref{fig:Duff_fd}b, or the time series data, see Figure~\ref{fig:Duff_fd}c, the trajectories of the nonlinear oscillator~\eqref{eq:Lin_fd} and the transformed trajectories of the linear system~\eqref{eq:duff_fd_lin} are practically indistinguishable. Moreover, the relative error remains below 6 percent (cf. Fig.~\ref{fig:Duff_fd}e).~\mycorr{ As mentioned in the previous section (cf. Section~\ref{sec:continuous}) the error  grows slowly due to the desynchronization of the true and approximated orbits, but generally remains bounded (cf. Appendix~\ref{app:error}).}

\mycorr{Additionally, we show the magnitudes and angles of two Koopman eigenfunctions of the forced-damped Duffing oscillator in Fig.~\ref{fig:Duff_fd}f. The magnitude of both eigenfunctions grows with the distance from the origin. Moreover, the phase of the first Koopman eigenfunction $\varphi_1$ associated with an eigenvalue with positive imaginary part increases in clockwise direction, while the phase of the second Koopman eigenfunction $\varphi_2$ decreases in clockwise direction. A similar behavior has been observed for the Koopman eigenfunctions of the nonlinear pendulum (cf. Fig.~\ref{fig:Pendulum}e.).}

\mycorr{Based on the observations from the preceding sections, we remark that introducing the coordinate $y_4$ in the linear system~\eqref{eq:Lin_fd} to distinguish the different stable steady state responses as well as the correct forcing frequency $\Omega$ are crucial to obtain accurate linearizations, while slight variation on the damping constant $c$ may be acceptable.}

The forced-damped Duffing equation~\eqref{eq:Duff_forced+damped} is an example of a nonlinear, single-degree-of-freedom oscillator. Systems with multiple degrees of freedom can also be studied with the techniques presented in this section. To this end, multiple degrees of freedom need to be added to the linear system~\eqref{eq:Lin_fd}. More specifically, for a $N$-degree-of-freedom system, an appropriate extension of the linear system~\eqref{eq:Lin_fd} would consist of the following $2N+2$ coordinates: $N$ coordinates for the positions, $N$ for the velocities, a phase coordinate for the periodic forcing and an additional coordinate if coexisting attractors are observed.

\section{Conclusions}

We have learned global linear embeddings of nonlinear systems using neural networks. Our approach shows that nonlinear phenomena can be captured by specific forms of the conjugated linear systems and appropriately parameterized transformations. We have exemplified this strategy on systems with continuous \mycorr{Koopman} spectra, such as the nonlinear pendulum and the Duffing oscillator, for which we have also learned a state-dependent frequency mapping. Furthermore, we have considered the Van-der-Pol oscillator, exhibiting a stable limit cycle, as well as the forced-damped Duffing oscillator with two coexisting steady states. Our analysis shows that both can be globally immersed into linear systems. For all systems the numerical computations show very good agreement of the nonlinear and the transformed, linear system.

While Koopman embeddings immerse nonlinear dynamics into an infinite dimensional space, our approach  yields a low-dimensional linearizing space. In the cases considered, the dimensions of the linear systems exceed the dimensions of the nonlinear system by one at most.   This low dimensionality leads to a considerable simplification in the learning of Koopman eigenfunctions and conjugate linear dynamics.

We have focused on three types of specific nonlinear behavior and illustrated our methodology on exemplary systems. \mycorr{The nonlinear phenomena  studied (continuous  \mycorr{Koopman} spectra, limit cycles, and coexisting steady states) occur in many applications, e.g.,  vortex induced vibrations~\cite{facchinetti2004coupling} (i.e. limit cycles, cf. Section~\ref{sec:LC}), nonlinear structural dynamics~\cite{kerschen2009nonlinear} (i.e. families of periodic orbits, cf. Section~\ref{sec:continuous}), or snap-through instabilities of MEMS-devices~\cite{das2009pull} (i.e. coexisting steady states, cf. Section~\ref{sec:FDDuffing}). Thus, it would be of interest to apply the proposed techniques to these applications.} Moreover, an extension to general systems with unknown limit sets may be envisioned. To this end, our approach advocates for phase space exploration and understanding first and then a subsequent linearization aided by a specifically selected linear system~\eqref{eq:lin_system}. Therein, the model systems analyzed could provide a starting point to collect a library of linear reference dynamics along with their immersions.

Moreover, the discussion from Section~\ref{sec:LC} can be extended to embedded quasi-periodic attractors in higher dimensional spaces. Similarly, the techniques from Sections~\ref{sec:continuous} and~\ref{sec:FDDuffing} could be extended to handle nested sequences of quasi-periodic tori or more complex excitation profiles. Future research could also include extensions to more complex attractor geometries. Especially, systems with internal resonances that can feature complex response patterns, such as, nested sequences of quasi-periodic orbits, quasi-periodic attractors, and chaos, are an area of future interest. It would be interesting to see if a low-dimensional version of a strange attractor could be learned along the same lines as presented here, or if chaotic systems require an infinite-dimensional representation, or if those systems are inherently non-linearizable.

%\todo[inline, color = green]{Maybe some more outlook, ideas? 1-2 points.}

\hspace{2cm}

\paragraph*{\textbf{Funding} } The authors declare that no funds, grants, or other support were received 
during the preparation of this manuscript.

\paragraph*{\textbf{Competing Interests}} The authors have no relevant financial or non-financial interests to disclose.

\paragraph*{\textbf{Author Contributions}}  T.B. and FK. reviewed the literature and conceptualized the approach. T.B. carried out the computations. T.B. and F.K. discussed the results and wrote the paper.

\paragraph*{\textbf{Data Availability}} The codes, generated data, and trained networks will be made publicly available.

\newpage

%\section{Nonlinear Pendulum}

%\begin{figure}  \centering    \includegraphics[scale=0.4]{IMG_7953.jpg}    \caption{Within the heteroclinic connection, the dynamics can be conjugated to the harmonic oscillator, outside, they are conjugated to a shear flow. }    \label{}\end{figure}There are two linear dynamics depending on the region:\begin{equation}\label{harmonic}   \begin{split}        \dot{x} & = -\omega y,\\        \dot{y} & = \omega x,    \end{split}\end{equation}within the heteroclinic connection and\begin{equation}\label{shear}\begin{split}    \dot{x} & = ax +  by,\\    \dot{y} & = 0.    \end{split}\end{equation}Since any linear similarity transform has to  preserve the spectrum, \eqref{harmonic} and \eqref{shear} are  distinct. For $a=0$, equation \eqref{shear} is also Hamiltonian, but is not associated to a fixed-point. So, the combined Hamiltonian would look something like\begin{equation}   \tilde{H}(x_1,y_1,x_2,y_2) = \frac{\omega}{2}(x_1^2+y_1^2)-\frac{b}{2}y_2^2\end{equation}for the standard symplectic form. The initial condition determines on which invariant hyperplane we evolve (and then non-linearly map to the true dynamics). This is that, outside the heteroclinic \textbf{Remark.}    A completely integrable system can, of course, always be completely linearized and the distinguished observables are the conserved quantities:    \begin{equation}        \frac{d}{dt} I_k = 0    \end{equation}

\bibliographystyle{abbrv}
\bibliography{Kooopman_bib}

\begin{thebibliography}{10}

\bibitem{abraham2019active}
I.~Abraham and T.~D. Murphey.
\newblock Active learning of dynamics for data-driven control using {K}oopman operators.
\newblock {\em IEEE Transactions on Robotics}, 35(5):1071--1083, 2019.

\bibitem{arbabi2017ergodic}
H.~Arbabi and I.~Mezic.
\newblock Ergodic theory, dynamic mode decomposition, and computation of spectral properties of the {K}oopman operator.
\newblock {\em SIAM Journal on Applied Dynamical Systems}, 16(4):2096--2126, 2017.

\bibitem{arnol2013mathematical}
V.~I. Arnol'd.
\newblock {\em Mathematical methods of classical mechanics}, volume~60.
\newblock Springer Science \& Business Media, 2013.

\bibitem{baladi2000positive}
V.~Baladi.
\newblock {\em Positive transfer operators and decay of correlations}, volume~16.
\newblock World scientific, 2000.

\bibitem{bauerlein2021phase}
B.~B{\"a}uerlein and K.~Avila.
\newblock Phase lag predicts nonlinear response maxima in liquid-sloshing experiments.
\newblock {\em Journal of Fluid Mechanics}, 925:A22, 2021.

\bibitem{bevanda2021koopman}
P.~Bevanda, S.~Sosnowski, and S.~Hirche.
\newblock Koopman operator dynamical models: Learning, analysis and control.
\newblock {\em Annual Reviews in Control}, 52:197--212, 2021.

\bibitem{brunton2016koopman}
S.~L. Brunton, B.~W. Brunton, J.~L. Proctor, and J.~N. Kutz.
\newblock Koopman invariant subspaces and finite linear representations of nonlinear dynamical systems for control.
\newblock {\em PloS one}, 11(2):e0150171, 2016.

\bibitem{brunton2021modern}
S.~L. Brunton, M.~Budi{\v{s}}i{\'c}, E.~Kaiser, and J.~N. Kutz.
\newblock Modern {K}oopman theory for dynamical systems.
\newblock {\em arXiv preprint arXiv:2102.12086}, 2021.

\bibitem{budivsic2012geometry}
M.~Budi{\v{s}}i{\'c} and I.~Mezi{\'c}.
\newblock Geometry of the ergodic quotient reveals coherent structures in flows.
\newblock {\em Physica D: Nonlinear Phenomena}, 241(15):1255--1269, 2012.

\bibitem{cenedese2022data}
M.~Cenedese, J.~Ax{\aa}s, B.~B{\"a}uerlein, K.~Avila, and G.~Haller.
\newblock Data-driven modeling and prediction of non-linearizable dynamics via spectral submanifolds.
\newblock {\em Nature communications}, 13(1):872, 2022.

\bibitem{das2009pull}
K.~Das and R.~Batra.
\newblock Pull-in and snap-through instabilities in transient deformations of microelectromechanical systems.
\newblock {\em Journal of Micromechanics and Microengineering}, 19(3):035008, 2009.

\bibitem{dormand1980family}
J.~R. Dormand and P.~J. Prince.
\newblock A family of embedded runge-kutta formulae.
\newblock {\em Journal of computational and applied mathematics}, 6(1):19--26, 1980.

\bibitem{eldering2018global}
J.~Eldering, M.~Kvalheim, and S.~Revzen.
\newblock Global linearization and fiber bundle structure of invariant manifolds.
\newblock {\em Nonlinearity}, 31(9):4202, 2018.

\bibitem{facchinetti2004coupling}
M.~L. Facchinetti, E.~De~Langre, and F.~Biolley.
\newblock Coupling of structure and wake oscillators in vortex-induced vibrations.
\newblock {\em Journal of Fluids and structures}, 19(2):123--140, 2004.

\bibitem{guckenheimer2013nonlinear}
J.~Guckenheimer and P.~Holmes.
\newblock {\em Nonlinear oscillations, dynamical systems, and bifurcations of vector fields}, volume~42.
\newblock Springer Science \& Business Media, 2013.

\bibitem{hislop2012introduction}
P.~D. Hislop and I.~M. Sigal.
\newblock {\em Introduction to spectral theory: With applications to {S}chr{\"o}dinger operators}, volume 113.
\newblock Springer Science \& Business Media, 2012.

\bibitem{holmes1997low}
P.~J. Holmes, J.~L. Lumley, G.~Berkooz, J.~C. Mattingly, and R.~W. Wittenberg.
\newblock Low-dimensional models of coherent structures in turbulence.
\newblock {\em Physics Reports}, 287(4):337--384, 1997.

\bibitem{hsu2013cell}
C.~S. Hsu.
\newblock {\em Cell-to-cell mapping: a method of global analysis for nonlinear systems}, volume~64.
\newblock Springer Science \& Business Media, 2013.

\bibitem{huang2018data}
B.~Huang and U.~Vaidya.
\newblock Data-driven approximation of transfer operators: Naturally structured dynamic mode decomposition.
\newblock In {\em 2018 Annual American Control Conference (ACC)}, pages 5659--5664. IEEE, 2018.

\bibitem{kaiser2021data}
E.~Kaiser, J.~N. Kutz, and S.~L. Brunton.
\newblock Data-driven discovery of {K}oopman eigenfunctions for control.
\newblock {\em Machine Learning: Science and Technology}, 2(3):035023, 2021.

\bibitem{kaplan2012understanding}
D.~Kaplan and L.~Glass.
\newblock {\em Understanding nonlinear dynamics}.
\newblock Springer Science \& Business Media, 2012.

\bibitem{katok1995introduction}
A.~Katok and B.~Hasselblatt.
\newblock {\em Introduction to the modern theory of dynamical systems}.
\newblock Cambridge University Press, 1995.

\bibitem{kelley1999iterative}
C.~T. Kelley.
\newblock {\em Iterative methods for optimization}.
\newblock SIAM, 1999.

\bibitem{kerschen2009nonlinear}
G.~Kerschen, M.~Peeters, J.-C. Golinval, and A.~F. Vakakis.
\newblock Nonlinear normal modes, part i: A useful framework for the structural dynamicist.
\newblock {\em Mechanical systems and signal processing}, 23(1):170--194, 2009.

\bibitem{kerschen2014modal}
G.~Kerschen, S.~Shaw, C.~Touz{\'e}, O.~Gendelman, B.~Cochelin, and A.~Vakakis.
\newblock {\em Modal analysis of nonlinear mechanical systems}, volume 555.
\newblock Springer, 2014.

\bibitem{klus2018data}
S.~Klus, F.~N{\"u}ske, P.~Koltai, H.~Wu, I.~Kevrekidis, C.~Sch{\"u}tte, and F.~No{\'e}.
\newblock Data-driven model reduction and transfer operator approximation.
\newblock {\em Journal of Nonlinear Science}, 28:985--1010, 2018.

\bibitem{koopman1931hamiltonian}
B.~O. Koopman.
\newblock Hamiltonian systems and transformation in {H}ilbert space.
\newblock {\em Proceedings of the National Academy of Sciences}, 17(5):315--318, 1931.

\bibitem{korda2018linear}
M.~Korda and I.~Mezi{\'c}.
\newblock Linear predictors for nonlinear dynamical systems: {K}oopman operator meets model predictive control.
\newblock {\em Automatica}, 93:149--160, 2018.

\bibitem{korda2020optimal}
M.~Korda and I.~Mezi{\'c}.
\newblock Optimal construction of {K}oopman eigenfunctions for prediction and control.
\newblock {\em IEEE Transactions on Automatic Control}, 65(12):5114--5129, 2020.

\bibitem{kovacic2011duffing}
I.~Kovacic and M.~J. Brennan.
\newblock {\em The {D}uffing equation: nonlinear oscillators and their behaviour}.
\newblock John Wiley \& Sons, 2011.

\bibitem{kowalski1991nonlinear}
K.~Kowalski and W.-H. Steeb.
\newblock {\em Nonlinear dynamical systems and {C}arleman linearization}.
\newblock World Scientific, 1991.

\bibitem{kvalheim2023linearizability}
M.~D. Kvalheim and P.~Arathoon.
\newblock Linearizability of flows by embeddings.
\newblock {\em arXiv preprint arXiv:2305.18288}, 2023.

\bibitem{lan2013linearization}
Y.~Lan and I.~Mezi{\'c}.
\newblock Linearization in the large of nonlinear systems and {K}oopman operator spectrum.
\newblock {\em Physica D: Nonlinear Phenomena}, 242(1):42--53, 2013.

\bibitem{liu2023properties}
Z.~Liu, N.~Ozay, and E.~D. Sontag.
\newblock Properties of immersions for systems with multiple limit sets with implications to learning {K}oopman embeddings.
\newblock {\em arXiv preprint arXiv:2312.17045}, 2023.

\bibitem{lorenz1963deterministic}
E.~N. Lorenz.
\newblock Deterministic nonperiodic flow.
\newblock {\em Journal of atmospheric sciences}, 20(2):130--141, 1963.

\bibitem{lusch2018deep}
B.~Lusch, J.~N. Kutz, and S.~L. Brunton.
\newblock Deep learning for universal linear embeddings of nonlinear dynamics.
\newblock {\em Nature communications}, 9(1):4950, 2018.

\bibitem{mauroy2020koopman}
A.~Mauroy, Y.~Susuki, and I.~Mezic.
\newblock {\em Koopman operator in systems and control}.
\newblock Springer, 2020.

\bibitem{mezic2020spectrum}
I.~Mezi{\'c}.
\newblock Spectrum of the {K}oopman operator, spectral expansions in functional spaces, and state-space geometry.
\newblock {\em Journal of Nonlinear Science}, 30(5):2091--2145, 2020.

\bibitem{murdock2003normal}
J.~A. Murdock.
\newblock {\em Normal forms and unfoldings for local dynamical systems}.
\newblock Springer, 2003.

\bibitem{nayfeh2008applied}
A.~H. Nayfeh and B.~Balachandran.
\newblock {\em Applied nonlinear dynamics: analytical, computational, and experimental methods}.
\newblock John Wiley \& Sons, 2008.

\bibitem{nayfeh2008nonlinear}
A.~H. Nayfeh and D.~T. Mook.
\newblock {\em Nonlinear oscillations}.
\newblock John Wiley \& Sons, 2008.

\bibitem{page2019koopman}
J.~Page and R.~R. Kerswell.
\newblock Koopman mode expansions between simple invariant solutions.
\newblock {\em Journal of Fluid Mechanics}, 879:1--27, 2019.

\bibitem{pinkus1999approximation}
A.~Pinkus.
\newblock Approximation theory of the {MLP} model in neural networks.
\newblock {\em Acta numerica}, 8:143--195, 1999.

\bibitem{rostamijavanani2023study}
A.~Rostamijavanani, S.~Li, and Y.~Yang.
\newblock A study on data-driven identification and representation of nonlinear dynamical systems with a physics-integrated deep learning approach: {K}oopman operators and nonlinear normal modes.
\newblock {\em Communications in Nonlinear Science and Numerical Simulation}, 123:107278, 2023.

\bibitem{ruelle2004thermodynamic}
D.~Ruelle.
\newblock {\em Thermodynamic formalism: the mathematical structure of equilibrium statistical mechanics}.
\newblock Cambridge University Press, 2004.

\bibitem{schmid2010dynamic}
P.~J. Schmid.
\newblock Dynamic mode decomposition of numerical and experimental data.
\newblock {\em Journal of fluid mechanics}, 656:5--28, 2010.

\bibitem{schmidhuber2015deep}
J.~Schmidhuber.
\newblock Deep learning in neural networks: An overview.
\newblock {\em Neural networks}, 61:85--117, 2015.

\bibitem{sijbrand1985properties}
J.~Sijbrand.
\newblock Properties of center manifolds.
\newblock {\em Transactions of the American Mathematical Society}, 289(2):431--469, 1985.

\bibitem{smale1967differentiable}
S.~Smale.
\newblock Differentiable dynamical systems.
\newblock {\em Bulletin of the American mathematical Society}, 73(6):747--817, 1967.

\bibitem{sternberg1958structure}
S.~Sternberg.
\newblock On the structure of local homeomorphisms of {E}uclidean n-space, ii.
\newblock {\em American Journal of Mathematics}, 80(3):623--631, 1958.

\bibitem{verhulst2006nonlinear}
F.~Verhulst.
\newblock {\em Nonlinear differential equations and dynamical systems}.
\newblock Springer Science \& Business Media, 2006.

\bibitem{whitham2011linear}
G.~B. Whitham.
\newblock {\em Linear and nonlinear waves}.
\newblock John Wiley \& Sons, 2011.

\bibitem{wiggins1994normally}
S.~Wiggins.
\newblock {\em Normally hyperbolic invariant manifolds in dynamical systems}, volume 105.
\newblock Springer Science \& Business Media, 1994.

\bibitem{williams2016extending}
M.~O. Williams, M.~S. Hemati, S.~T. Dawson, I.~G. Kevrekidis, and C.~W. Rowley.
\newblock Extending data-driven {K}oopman analysis to actuated systems.
\newblock {\em IFAC-PapersOnLine}, 49(18):704--709, 2016.

\bibitem{williams2015data}
M.~O. Williams, I.~G. Kevrekidis, and C.~W. Rowley.
\newblock A data--driven approximation of the {K}oopman operator: Extending dynamic mode decomposition.
\newblock {\em Journal of Nonlinear Science}, 25:1307--1346, 2015.

\bibitem{williams2014kernel}
M.~O. Williams, C.~W. Rowley, and I.~G. Kevrekidis.
\newblock A kernel-based approach to data-driven {K}oopman spectral analysis.
\newblock {\em arXiv preprint arXiv:1411.2260}, 2014.

\bibitem{yeung2019learning}
E.~Yeung, S.~Kundu, and N.~Hodas.
\newblock Learning deep neural network representations for {K}oopman operators of nonlinear dynamical systems.
\newblock In {\em 2019 American Control Conference (ACC)}, pages 4832--4839. IEEE, 2019.

\end{thebibliography}

\clearpage

\appendix 

\section{\mycorr{Upper bound on the Approximation Error}}
\label{app:error}
\mycorr{Approximating the mappings $\boldsymbol{\Phi}$ and $\boldsymbol{\Phi}^{-1}$ with neural networks (cf.~Section\ref{sec:learning} introduces differences between the trajectories of the nonlinear system $\mathbf{x}$ and the transformed linear system $\bm{\Phi}^{-1}(\mathbf{y}(t))=\bm{\Phi}^{-1}(e^{\mathbf{A}t}\bm{\Phi}(\mathbf{x}_0))$ (cf. equation~\eqref{eq:x_lin}). We denote the approximations with $\tilde{\boldsymbol{\Phi}}$ and $\tilde{\boldsymbol{\Phi}}^{-1}$ and obtain for the difference between the trajectories $\mathbf{x}(t)$ and $\bm{\Phi}^{-1}(\mathbf{y}(t))$ 
   \begin{widetext}
   \begin{equation}
   	\label{eq:error_bound}
  	\begin{split}
 	\varepsilon(t)&:=|\mathbf{x}(t)-\tilde{\boldsymbol{\Phi}}^{-1}(\mathbf{F}^t(\tilde{\boldsymbol{\Phi}}(\mathbf{x}_0)))| \\
 	&\leq |\mathbf{x}(t)-\boldsymbol{\Phi}^{-1}(\mathbf{F}^t(\boldsymbol{\Phi}(\mathbf{x}_0)))|
 	+|\tilde{\boldsymbol{\Phi}}^{-1}(\mathbf{F}^t(\tilde{\boldsymbol{\Phi}}(\mathbf{x}_0)))-\boldsymbol{\Phi}^{-1}(\mathbf{F}^t(\boldsymbol{\Phi}(\mathbf{x}_0)))|
 	\\
 	&\leq|\tilde{\boldsymbol{\Phi}}^{-1}(\mathbf{F}^t(\boldsymbol{\Phi}(\mathbf{x}_0)))-\boldsymbol{\Phi}^{-1}(\mathbf{F}^t(\boldsymbol{\Phi}(\mathbf{x}_0)))|+ |\tilde{\boldsymbol{\Phi}}^{-1}(\mathbf{F}^t(\tilde{\boldsymbol{\Phi}}(\mathbf{x}_0)))-\tilde{\boldsymbol{\Phi}}^{-1}(\mathbf{F}^t(\boldsymbol{\Phi}(\mathbf{x}_0))) |\\
 	&\leq |\tilde{\boldsymbol{\Phi}}^{-1}(  \mathbf{y}(t))-\boldsymbol{\Phi}^{-1}(\mathbf{y}(t))|+C_{\tilde{\boldsymbol{\Phi}}^{-1}}|\mathbf{F}^t(\tilde{\boldsymbol{\Phi}}(\mathbf{x}_0))-\mathbf{F}^t(\boldsymbol{\Phi}(\mathbf{x}_0))|\\
 	&\leq |\tilde{\boldsymbol{\Phi}}^{-1}(  \mathbf{y}(t))-\mathbf{x}(t)|+C_{\tilde{\boldsymbol{\Phi}}^{-1}}|e^{\mathbf{A}t}(\tilde{\boldsymbol{\Phi}}(\mathbf{x}_0)-\boldsymbol{\Phi}(\mathbf{x}_0))| \\
 	&\leq |\tilde{\boldsymbol{\Phi}}^{-1}(  \mathbf{y})- \mathbf{x}|+C_{\tilde{\boldsymbol{\Phi}}^{-1}} \,||e^{\mathbf{A}t}|| \, |\tilde{\boldsymbol{\Phi}}(\mathbf{x}_0)-\boldsymbol{\Phi}(\mathbf{x}_0))| \\
 	&\leq \underbrace{|\tilde{\boldsymbol{\Phi}}^{-1}(\mathbf{y})- \mathbf{x}|}_{\text{Residual from~\eqref{eq:minimizations}}}+C_{\tilde{\boldsymbol{\Phi}}^{-1}} \,||e^{\mathbf{A}t}|| \, \underbrace{|\tilde{\boldsymbol{\Phi}}(\mathbf{x}_0)-\mathbf{y}_0)|}_{\text{Residual from~\eqref{eq:minimizations}}}, 
 	\end{split}
 \end{equation} 
 \end{widetext}
 where $C_{\tilde{\boldsymbol{\Phi}}^{-1}}$ denotes the Lipschitz constant of the neural network $\tilde{\boldsymbol{\Phi}}^{-1}$. The deep neural network $\tilde{\boldsymbol{\Phi}}^{-1}$ is Lipschitz continuous, since $\tilde{\boldsymbol{\Phi}}^{-1}$ is a  composition of infinitely often differentiable functions (i.e. linear transformations and $\mbox{tansig}(x)$ functions). From the last line of equation~\eqref{eq:error_bound} we deduce that the error generally grows with the linear dynamics~\eqref{eq:lin_system}. Moreover, the upper bound~\eqref{eq:error_bound} scales linearly with the residuals of the minimizations~\eqref{eq:minimizations}. }

\section{Examples of Explicit Immersions for Some Nonlinear Systems}
\label{app:examples}

In this appendix, we collect some explicit examples of transformations and embeddings from nonlinear model systems into linear systems.

\subsection{Continuous \mycorr{Koopman} Spectrum: Transform to the Linear Harmonic Oscillator}
\label{app:harmonic}

For systems with families of nested periodic orbits, i.e., continuous spectra of the Koopman operator, we start with the linear dynamics in observables and construct nonlinearly equivalent systems. To this end, consider the linear harmonic oscillator 
\begin{equation}\label{harmonicosc}
\begin{split}
    \dot{y}_1&=\omega y_2,\\
    \dot{y}_2&=-\omega y_1,\\
\end{split}
\end{equation}
with frequency $\omega>0$ in the space of observables. We write
\begin{equation}
    \begin{split}
        r & =\sqrt{y_1^2+y_2^2},\\
        \theta &  =\arctan \left( y_1/y_2\right),
    \end{split}
\end{equation}
for the polar coordinates in observable space. The inverse embedding 
\begin{equation}
\begin{split}
x_1&=\psi_1(y_1,y_2)=f_r\left(r,\theta\right),\\
x_2&=\psi_2(y_1,y_2)=f_{\phi}\left(r,\theta\right),\\
\end{split}
\end{equation}
for any two functions $\psi_1, \psi_2$ and $f_r,f_{\phi}$, respectively, leads to the nonlinear planar system 
\begin{equation}
\label{continuousmodel}
\begin{split}
\dot{x}_1&=\partial_{y_1}\psi_1(y_1,y_2)\omega y_2-\partial_{y_2}\psi_1(y_1,y_2)\omega y_1,\\
\dot{x}_2&=\partial_{y_1}\psi_2(y_1,y_2)\omega y_2-\partial_{y_2}\psi_2(y_1,y_2)\omega y_1.\\
\end{split}
\end{equation}
Any system of the form~\eqref{continuousmodel} can hence be transformed to the linear harmonic oscillator~\eqref{harmonicosc}. Indeed, system ~\eqref{continuousmodel} yields families of periodic orbits with an instantaneous angular speed given by $f_{\phi}(r,\theta)$ and radius  $f_{r}(r,\theta)$. 

The choice $f_r=r$ and $f_{\phi}(r,\theta)=\theta-\theta\sin(\pi/2\cdot r)$, for example, yields a systems with radially varying frequency $\Omega=\omega-\omega\sin(\pi/2\cdot r)$. For $r=0$, the frequency is equal to the linear frequency $\omega$ and decreases with increasing $r$ until it reaches zero at $r=1$.

\subsection{Limit Cycle}
\label{app:LC}
For the nonlinear planar system
\begin{equation}
\label{eq:LC}
\begin{split}
    \dot{x}_1&=\frac{x_1}{\sqrt{x_1^2+x_2^2}}+(x_2-x_1),\\
    \dot{x}_2&=\frac{x_2}{\sqrt{x_1^2+x_2^2}}-(x_1+x_2),
\end{split}
\end{equation}
consider the radius and the angle as observables,
 \begin{equation}
\label{eq:LC_obs}
\begin{split}
y_1&=\varphi_1(x_1,x_2)=\sqrt{x_1^2+x_2^2}\\
y_2&=\varphi_2(x_1,x_2)=\arctan\left(\frac{x_1}{x_2}\right).
\end{split}
\end{equation}
For the first observable $y_1$ we calculate
\begin{equation}
\label{eq:LC_X1_dyn}
\begin{split}
\dot{y}_1&=\frac{\dot{x}_1x_1+\dot{x}_2x_2}{\sqrt{x_1^2+x_2^2}} \\
                    &=\frac{\frac{x_1^2}{\sqrt{x_1^2+x_2^2}}+(x_2-x_1)x_1+\frac{x_2^2}{\sqrt{x_1^2+x_2^2}}-(x_1+x_2)x_2}{\sqrt{x_1^2+x_2^2}} \\
                    &=\frac{\frac{x_1^2+x_2^2}{\sqrt{x_1^2+x_2^2}}-x_1^2-x_2^2}{\sqrt{x_1^2+x_2^2}} \\
                    &=\frac{\frac{y_1^2}{y_1}-y_1^2}{y_1}=-y_1+1 
\end{split}
\end{equation}
as well as for the second observable $y_2$,
\begin{small}
\begin{equation}
\label{eq:LC_X2_dyn}
\begin{split}
\dot{y}_2&=\frac{1}{1+\frac{x_1^2}{x_2^2}}\frac{\dot{x}_1x_2-x_1\dot{x}_2}{x_2^2}\\
                    &= \frac{1}{1+\frac{x_1^2}{x_2^2}}\frac{\frac{x_1x_2}{\sqrt{x_1^2+x_2^2}}+(x_2-x_1)x_2-\frac{x_2x_1}{\sqrt{x_1^2+x_2^2}}+(x_1+x_2)x_1}{x_2^2}\\
                    &= \frac{(x_2^2+x_1^2)}{x_1^2+x_2^2}=1. 
\end{split}
\end{equation}
\end{small}
Writing system~\eqref{eq:LC} in polar coordinates with the observables~\eqref{eq:LC_obs} yields the dynamics
\begin{equation}
\label{eq:LC_koopman}
\begin{split}
    \dot{y}_1&=-y_1+1,\\
    \dot{y}_2&=1.\\
\end{split}
\end{equation}
From equation~\eqref{eq:LC_koopman}, we see immediately that system ~\eqref{eq:LC_obs} has the unit circle as a limit cycle. Introducing another dummy observable $\dot{y_3} = 0$ gives the three-dimensional linearization 
\begin{equation}
   \begin{split}
    \dot{y}_1&=-y_1+y_3,\\
    \dot{y}_2&=y_3,\\
    \dot{y}_3&=0.
\end{split} 
\end{equation}
For the special solution $y_3=1$, we recover system~\eqref{eq:LC}.

\subsection{Coexisting steady state}
\label{app:coexist}
In this section, we provide a simple example of a planar nonlinear system with multiple fixed-points that can be embedded into a linear system. Indeed, consider the planar nonlinear system 
\begin{subequations}
\label{eq:coexsist}
\begin{align}
\label{eq:coexsist_rad}
   \dot{x}_1&=x_1(\lambda_1 \sin(x_2)^2+\lambda_2\cos(x_1)^2), \\
   \label{eq:coexsist_phi}
    \dot{x}_2&=(\lambda_1-\lambda_2)\sin(x_2)\cos(x_2) ,
\end{align}
\end{subequations}
for $\lambda_2<\lambda_1<0$, which has the fixed points
\begin{equation}
    \begin{split}
    x_1&=0,\qquad x_2=0,\qquad ~~~~\Rightarrow ~~\text{stable},\\
    x_1&=0,\qquad x_2=\pi,\qquad ~~~~\Rightarrow ~~\text{stable}, \\
    x_1&=0,\qquad x_2=\pi/2,\qquad ~\Rightarrow ~~\text{unstable},\\
    x_1&=0,\qquad x_2=3\pi/2,\qquad \Rightarrow ~~\text{unstable}.
    \end{split}
\end{equation}
We define the observables 
\begin{equation}
\label{eq:coex_obs}
    \begin{split}
        y_1&=x_1\sin(x_2)\\
        y_2&=x_1\cos(x_2)
    \end{split}
\end{equation}
Calculating the time derivative of the first observable yields
\begin{equation}
\label{eq:coex_y1_dyn}
\begin{split}
\dot{y}_1&=\dot{x}_1\sin(x_2)+x_1\dot{x}_2\cos(x_2)\\
                    &=x_1(\lambda_1 \sin(x_2)^2+\lambda_2\cos(x_1)^2)\sin(x_2)\\
                    &\quad  +x_1(\lambda_1-\lambda_2)\sin(x_2)\cos(x_2)\cos(x_2)\\
                    &=x_1\lambda_1 \sin(x_2)^3+x_1\lambda_1\sin(x_2)\cos(x_2)^2\\
                    &=x_1\lambda_1 \sin(x_2)=\lambda_1y_1.\\
\end{split}
\end{equation}
For the second observable, we have 
\begin{equation}
\label{eq:coex_X2_dyn}
\begin{split}
\dot{y}_2&=\dot{x}_1\cos(x_2)-x_1\dot{x}_2\sin(x_2)\\
                    &=x_1(\lambda_1 \sin(x_2)^2+\lambda_2\cos(x_1)^2)\cos(x_2)\\
                    &\quad -x_1(\lambda_1-\lambda_2)\sin(x_2)\cos(x_2)\sin(x_2)\\
                    &=x_1\lambda_2 \cos(x_2)^3+x_1\lambda_2\sin(x_2)^2\cos(x_2)\\
                    &=x_1\lambda_2 \cos(x_2)=\lambda_2y_2,.\\
\end{split}
\end{equation}
Consequently, the observables~\eqref{eq:coex_obs} satisfy the linear dynamics 
\begin{equation}
    \label{eq:lin_FP}
    \begin{split}
    \dot{y}_1&=\lambda_1y_1,\\
    \dot{y}_2&=\lambda_2y_2. 
    \end{split}
\end{equation}

Indeed, the nonlinear system~\eqref{eq:coexsist_rad} is obtained by transforming the linear system ~\eqref{eq:lin_FP} into polar coordinates.

\begin{figure}[ht!]
        \begin{center}
           \includegraphics[width=0.45\textwidth]{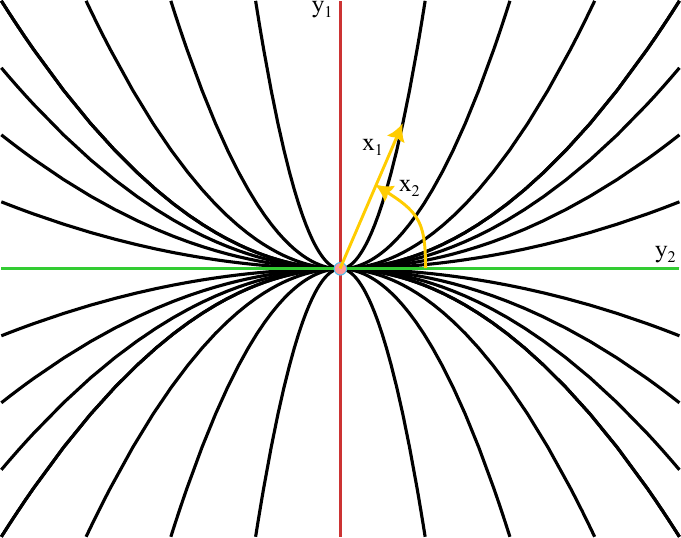}
 	\end{center}
 	\caption{The two-dimensional linear system ~\eqref{eq:lin_FP} has two stable and two unstable fixed points in polar coordinates.}
 	\label{fig:2D_lin_sys}
 \end{figure}

Let us take a closer look at the coexistence of steady states. Figure ~\ref{fig:2D_lin_sys} shows the phase portrait of system ~\eqref{eq:lin_FP}.   The only fixed point of system~\eqref{eq:lin_FP} is the trivial solution $y_1=y_2=0$, which is stable since $\lambda_2<\lambda_1<0$. General trajectories of the linear system ~\eqref{eq:lin_FP} in polar coordinates (shown in black in Fig.~\ref{fig:2D_lin_sys}) decay along the slow stable direction towards the origin. If the initial condition for the coordinate $y_1$ is positive, then trajectories approach the fixed point along the line $\varphi=0$. Thus, in polar coordinates these trajectories approach a steady state with an angle $\varphi=x_1=0$. If, on the other hand, the initial condition for the coordinate $y_1$ is negative, trajectories approach the fixed point along the line $\varphi=x_1=\pi$.

For the initial conditions $y_1=0$ and  $y_2>0$, trajectories decay along the fast stable direction. In polar coordinates, this manifold corresponds to the line $\varphi=x_1=\pi/2$. This solution is unstable, as small perturbations causes trajectories to approach the fixed point along the slow stable direction. Analogously, the trajectory associated to the initial conditions $y_1=0$ and  $y_2<0$ approach a steady state angle $\varphi=x_1=3\pi/2$, which is also unstable. 

\section{\mycorr{Visualizations of the mappings $\bm{\Phi}$ and $\bm{\Phi}^{-1}$ for the nonlinear pendulum~\eqref{eq:pendulum}}}
\label{app:pendulum_phi_visu}
\begin{figure*}
 	\begin{center}
			\includegraphics[width=\textwidth]{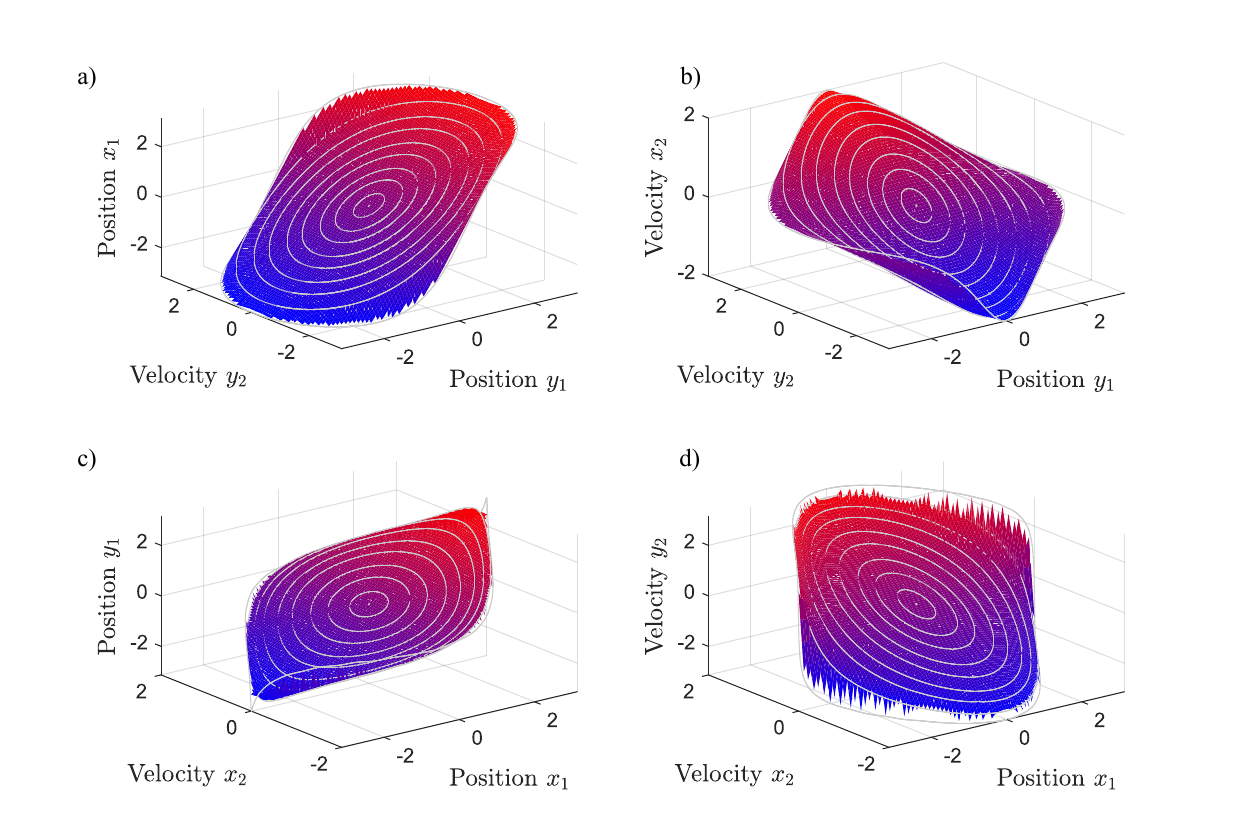}
 	\end{center}
 	\caption{Visualizations of the outputs of the deep neural networks. \textbf{a)} and \textbf{b)} $\bm{\Phi}^{-1}(\mathbf{y})$ maps coordinates of the linear oscillator~\eqref{eq:lin_osci} to coordinates of the nonlinear pendulum~\eqref{eq:pendulum}.  \textbf{c)} and \textbf{d)} $\bm{\Phi}(\mathbf{x})$ maps coordinates of the nonlinear pendulum~\eqref{eq:pendulum} to  coordinates of the linear oscillator~\eqref{eq:lin_osci}. Gray lines are trajectories of the training data. }
 	\label{fig:Pendulum_Phi_visu}
\end{figure*}

\mycorr{In Fig.~\ref{fig:Pendulum_Phi_visu} we visualize the output of the mappings $\bm{\Phi}$ and $\bm{\Phi}^{-1}$. More specifically, Fig.~\ref{fig:Pendulum_Phi_visu}a and Fig.~\ref{fig:Pendulum_Phi_visu}b show the transformation of coordinates of the linear oscillator~\eqref{eq:lin_osci} to coordinates of the nonlinear pendulum~\eqref{eq:pendulum}, while Fig.~\ref{fig:Pendulum_Phi_visu}c and Fig.~\ref{fig:Pendulum_Phi_visu}d show the mapping of  coordinates of the nonlinear pendulum~\eqref{eq:pendulum} to coordinates of the linear oscillator~\eqref{eq:lin_osci}.}

\section{Necessity of including the nonlinear period in the transformed, linear oscillator~\eqref{eq:lin_pendulum}} %
\label{app:Tcorr}
\begin{figure*}
 	\begin{center}
			\includegraphics[width=\textwidth]{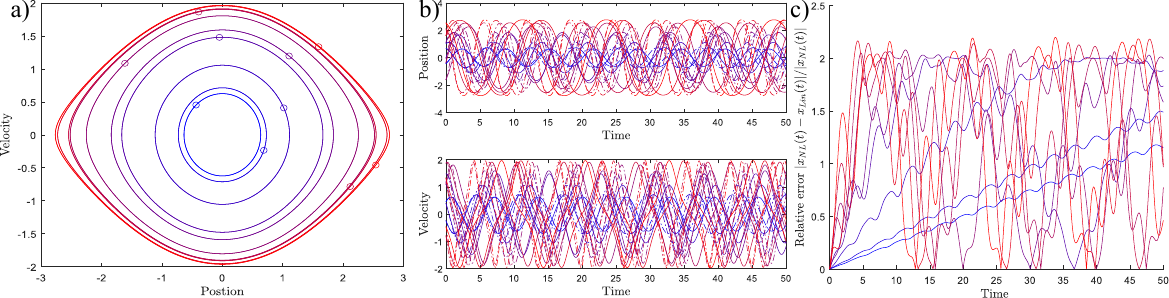}
 	\end{center}
 	\caption{Transformed linear oscillator without adjusting for the nonlinear period~\eqref{eq:lin_pendulum_noTcorr}. Solid lines demarcate the nonlinear pendulum~\eqref{eq:pendulum} and dashed-dotted lines the trajectories~\eqref{eq:lin_pendulum_noTcorr}.   \textbf{a)} Phase space. Circle markers indicate the randomly selected initial conditions.  \textbf{b)} Raw time series of the position and the velocities. \textbf{c)} Relative error between trajectories of the transformed linear oscillator and the nonlinear pendulum. }
 	\label{fig:Pendulum_noTcorr}
\end{figure*}

To emphasize the importance of adjusting the period in equation~\eqref{eq:lin_pendulum}, we omit this adjustment and plot the trajectories
\begin{equation}
\label{eq:lin_pendulum_noTcorr}
     \begin{split}
         (\tilde{r}_{\rm Lin}(t), \tilde{\theta}_{\rm Lin}(t)) =\bm{\Phi}^{-1}\left(r_l(t),\theta_l(t)\right),
     \end{split}
\end{equation}
in Fig.~\ref{fig:Pendulum_noTcorr}.  While the phase space geometry is accurately reconstructed in Fig.~\ref{fig:Pendulum_noTcorr}a, the raw time series generated via equation~\eqref{eq:lin_pendulum_noTcorr} deviate quickly from the solutions of the nonlinear pendulum~\eqref{eq:pendulum}. This shows that including the nonlinear period in equation~\eqref{eq:lin_pendulum} is indeed crucial.

\section{Relationship of the Koopman operator to linear immersions}
\label{app:Koopman}
The Koopman operator, acting on scalar observables $\psi:\mathbb{R}^n\to\mathbb{C}$, is defined as the semi-group of compositions $\mathcal{K}_t[\psi] = \psi\circ\mathbf{F}^t$ and its infinitesimal generator is the Lie-derivative along the right-hand side of the nonlinear system~\eqref{eq:NL_system}:
\begin{equation}\label{defL}
    \mathcal{L}\psi = \lim_{t\to 0+} \frac{1}{t} (\mathcal{K}_t\psi-\psi)=\nabla \psi\cdot \mathbf{f}.
\end{equation}
The eigenfunctions of the linear operator~\eqref{defL} satisfy $\mathcal{L}\psi = \lambda \psi$, for $\lambda\in\mathbb{C}$, and imply a corresponding one-parameter family of eigenvalues of the Koopman operator, $ \mathcal{K}_t\psi = e^{t\lambda }\psi$. Clearly, the product of two Koopman eigenfunctions defines an eigenfunction as well.

The Koopman theory is inherently linked to the linear conjugacy of flow maps and finite-dimensional immersions as follows. First, taking a time-derivative of equation~\eqref{immserion} and evaluating at $t=0$ leads to the relation
\begin{equation}\label{immersiondiff}
   \frac{\partial\bm{\Phi}}{\partial\bm{x}}(\bm{x}_0) \bm{f}(\bm{x}_0) = \bm{A} \bm{\Phi}(\bm{x}_0).
\end{equation}
Furthermore, we introduce the Jordan decomposition $\bm{\Lambda}$ of $\bm{A}$ as $\bm{A}=\bm{V}\bm{\Lambda}\bm{V}^{-1}$, where $\bm{V}$ denote the generalized eigenspaces. Let $\lambda_1,...,\lambda_M$ denote the eigenvalues of the matrix $\bm{A}$ (counted with multiplicity)  
and define coordinates $\bm{z}=\bm{V}^{-1}\bm{y}$, since for $\bm{z}=\bm{V}^{-1}\bm{\Phi}(\bm{x})$,
\begin{small}
\begin{equation}\label{yeigen}
    \begin{split}
        (\bm{f}(\bm{x})\cdot\nabla_{\bm{x}}) \bm{V}^{-1}\bm{\Phi}(\bm{x})
        %= (\bm{f}(\bm{x})\cdot\nabla_{\bm{x}}) \bm{\Phi}(\bm{x})\\
        =\bm{V}^{-1}\frac{\partial\bm{\Phi}}{\partial\bm{x}}\bm{f}(\bm{x})= \bm{\Lambda} \bm{V}^{-1}\bm{\Phi}(\bm{x}).
       % & = \bm{\Lambda} \bm{\Phi}(\bm{x})\\
    \end{split}
\end{equation}
\end{small}
Consequently, any finite-dimensional linear embedding defines a specific set of Koopman eigenfunctions through its coordinates.

\section{\mycorr{Visualizations of the mappings $\bm{\Phi}$ and $\bm{\Phi}^{-1}$ for the Van-der-Pol oscillator~\eqref{eq:VdP}}}
\label{app:VdP_phi_visu}
\begin{figure*}
 	\begin{center}
			\includegraphics[width=\textwidth]{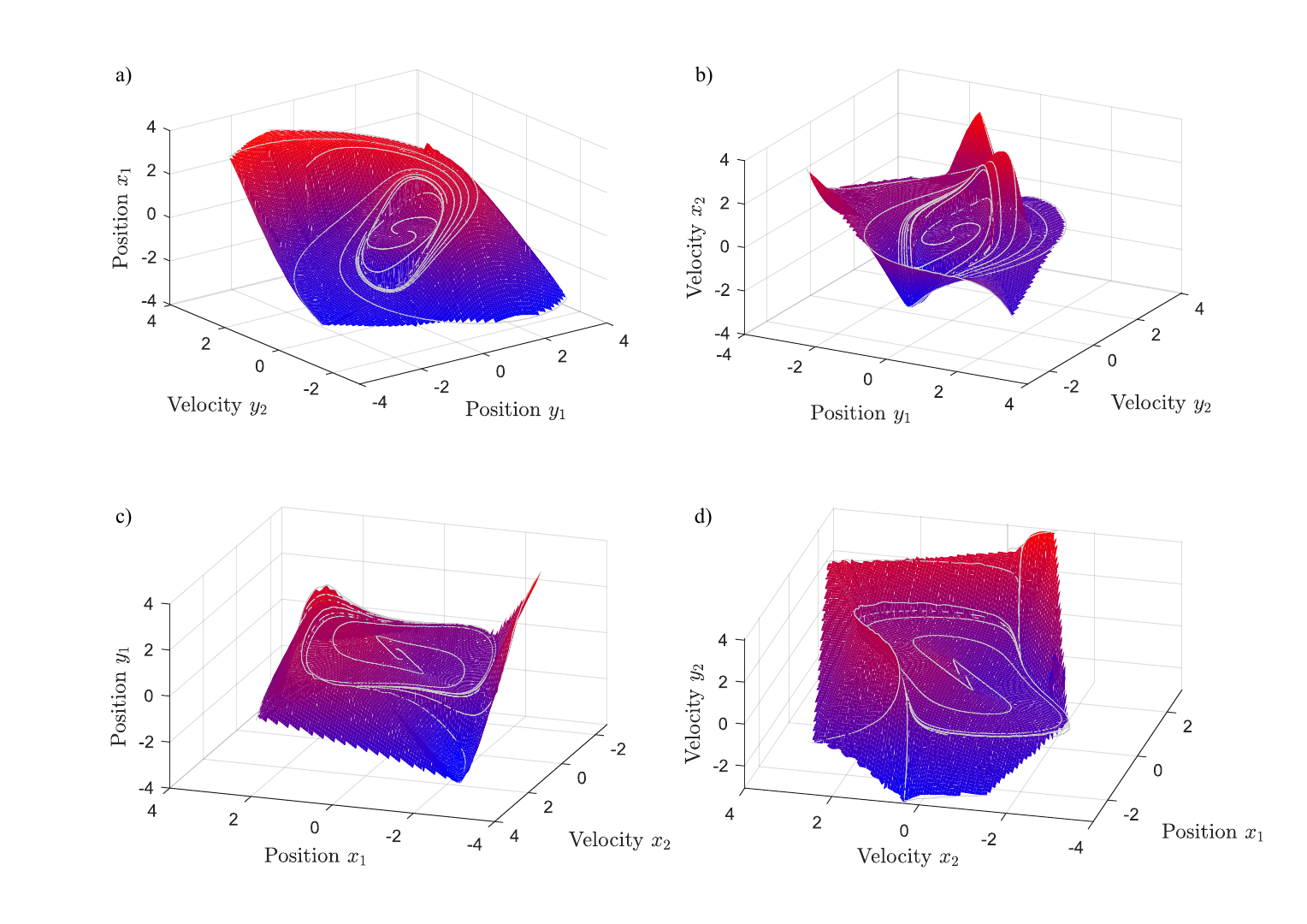}
 	\end{center}
 	\caption{Visualizations of the outputs of the deep neural networks. \textbf{a)} and \textbf{b)} $\bm{\Phi}^{-1}(\mathbf{y})$ maps coordinates of the linear system~\eqref{eq:lin_VdP} to coordinates of the Van-der-Pol oscillator~\eqref{eq:VdP}.  \textbf{c)} and \textbf{d)} $\bm{\Phi}(\mathbf{x})$ maps coordinates of the Van-der-Pol oscillator~\eqref{eq:VdP} to  coordinates of the linear system~\eqref{eq:lin_VdP}. Gray lines are trajectories of the training data.}
 	\label{fig:VdP_Phi_visu}
\end{figure*}

\mycorr{In Fig.~\ref{fig:VdP_Phi_visu} we visualize the output of the mappings $\bm{\Phi}$ and $\bm{\Phi}^{-1}$. More specifically, Fig.~\ref{fig:VdP_Phi_visu}a and Fig.~\ref{fig:VdP_Phi_visu}b show the transformation of coordinates of the linear oscillator~\eqref{eq:lin_VdP} to coordinates of the Van-der-Pol oscillator~\eqref{eq:VdP}, while Fig.~\ref{fig:VdP_Phi_visu}c and Fig.~\ref{fig:VdP_Phi_visu}d show the mapping of  coordinates of the Van-der-Pol oscillator~\eqref{eq:VdP} to coordinates of the linear system~\eqref{eq:lin_VdP}.}

\end{document}